\documentclass{article}
\usepackage{arxiv}
\usepackage[utf8]{inputenc} 
\usepackage[T1]{fontenc}    
\usepackage{hyperref}       
\usepackage{url}            
\usepackage{booktabs}       
\usepackage{amsfonts}       
\usepackage{nicefrac}       
\usepackage{microtype}      
\usepackage{graphicx}
\usepackage{doi}
\usepackage{graphicx}
\usepackage{mathptmx}      
\usepackage{amsmath}
\usepackage{amsthm}
\usepackage{amssymb}
\usepackage{xcolor}
\usepackage{bm}
\usepackage{enumerate}
\usepackage{latexsym}
\usepackage[ruled]{algorithm2e}
\usepackage{adjustbox}
\usepackage{hyperref}
\usepackage{url}
\usepackage{multirow}
\usepackage{soul}
\usepackage[square,numbers]{natbib}
\newtheorem{theorem}{Theorem}[section]

\newtheorem{lemma}[theorem]{Lemma}

\SetKwInput{KwInput}{Input}        
\SetKwInput{KwOutput}{Output}      

\usepackage{pifont}
\newcommand{\cmark}{\ding{51}}%
\newcommand{\xmark}{\ding{55}}
\newcommand{\erre}{\mathbb{R}}

\makeatletter
\newcommand{\pushright}[1]{\ifmeasuring@#1\else\omit\hfill$\displaystyle#1$\fi\ignorespaces}
\newcommand{\pushleft}[1]{\ifmeasuring@#1\else\omit$\displaystyle#1$\hfill\fi\ignorespaces}
\makeatother
\newcommand{\rev}[1]{\textcolor{black}{#1}}

\title{On the generation of Metric TSP instances with a large integrality gap by branch-and-cut}

\author{ Eleonora Vercesi\\
	Department of Mathematics\\
	University of Pavia\\
	\texttt{eleonora.vercesi01@universitadipavia.it} \\
	\And
	Stefano Gualandi\\
	Department of Mathematics\\
	University of Pavia\\
	\texttt{stefano.gualandi@unipv.it} \\
	\And
	Monaldo Mastrolilli\\
	IDISA, Lugano\\
	\texttt{monaldo.mastrolilli@idsia.ch} \\
    \And
	Luca Maria Gambardella\\
	Faculty of Informatics, USI, IDSIA - USI-SUPSI, Lugano\\
	\texttt{luca.gambardella@idsia.ch} \\

}

\renewcommand{\shorttitle}{On the generation of Metric TSP instances with a large integrality gap by branch-and-cut}

\hypersetup{
pdftitle={A template for the arxiv style},
pdfsubject={q-bio.NC, q-bio.QM},
pdfauthor={David S.~Hippocampus, Elias D.~Striatum},
pdfkeywords={First keyword, Second keyword, More},
}

\begin{document}
\maketitle

\begin{abstract}
This paper introduces a computational method for generating metric Travelling Salesman Problem (TSP) instances having a large integrality gap.
The method is based on the solution of an integer programming problem, called IH-OPT, that takes as input a fractional solution of the Subtour Elimination Problem (SEP) on a TSP instance and computes a TSP instance having an integrality gap larger than or equal to the integrality gap of the first instance.
The decision variables of IH-OPT are the entries of the TSP cost matrix, and the constraints are defined by the intersection of the metric cone with an exponential number of inequalities, one for each possible TSP tour.
Given the very large number of constraints, we have implemented a branch-and-cut algorithm for solving IH-OPT.
Then, by sampling cost vectors over the metric polytope and by solving the corresponding SEP, we can generate random fractional vertices of the SEP polytope.
If we solve the IH-OPT problem for every sampled vertex using our branch-and-cut algorithm, we can select the generated TSP instance (i.e., cost vector), yielding the longest runtime for Concorde, the state-of-the-art TSP solver.
Our computational results show that our method is very effective in producing challenging instances.
As a by-product, we release the {\tt Hard-TSPLIB}, a library of 41 small metric TSP instances which have a large integrality gap and are challenging in terms of runtime for Concorde.

\end{abstract}

\section{Introduction}
\label{intro}
The Branch-and-Cut (B\&C) algorithms are powerful tools to solve NP-hard problems.
A core component of these algorithms is the solution of several Linear Programming (LP) relaxations that appear while searching for an integral optimal solution.
Given an instance of an ILP model, its {\it integrality gap} for one linear relaxation is the ratio between the optimal solution of the integer problem and the corresponding LP relaxation.
In practice, this is a measure of how much the relaxation is good to approximate the original problem.
As a rule of thumb, given an Integer Linear Programming (ILP) model for an NP-hard problem, the larger is the integrality gap of the initial LP relaxation, the longer is the runtime for a B\&C algorithm to prove that the best solution found is optimal.
In practice, a large integrality gap at the root node very often implies that a great number of branch-and-bound nodes must be visited before proving the optimality of an integral solution.
For this reason, a significant effort in designing efficient B\&C algorithms is spent looking for tight LP relaxations \cite{mitchell2002branch}.

Let us consider, for instance, the Travelling Salesman Problem (TSP):
Given a set of $n$ nodes $V=\{1,\dots,n\}$, a square matrix $\bm{C} \in \mathbb{R}^{n \times n}_+$ whose entries $c_{ij}$ represent the cost of going from node $i$ to node $j$, we have to find the cyclic permutation $\pi$ of $V$ such that the total cost $\sum_{i \in V} c_{i,\pi_i}$ is minimum.
If the matrix $\bm{C}$ is (i) {\it symmetric}, that is, $c_{ij}=c_{ji}$ for all $i,j \in V$, (ii) satisfies the {\it triangle inequalities} $c_{ij} + c_{jk} \geq c_{ik}, \forall i,j,k \in V$, and (iii) $c_{ij}=0$ if and only if $i=j$, then the matrix $\bm{C}$ defines a {\it metric} on $V$ and gives a {\it Metric TSP} instance.
Indeed, there exist several ILP models to solve the TSP \cite{applegate2006traveling,dantzig1954solution,miller1960integer,orman2007survey}.
However, the most successful exact method is the branch-and-cut algorithm introduced by Padberg and Rinaldi \cite{padberg1991branch}, which includes as a core component the solution of an LP relaxation based on the subtour elimination constraints \cite{padberg1990facet}.
This LP relaxation is known in the literature as the {\it Subtour Elimination Problem (SEP)}.
Nowadays, the state-of-the-art exact software for the symmetric TSP is Concorde \cite{applegate1998solution}, a solver which implements a sophisticated B\&C algorithm and which includes the separation of several facet-defining inequalities \cite{applegate2006traveling}.
Concorde has a track record in solving the instances from the {\tt TSPLIB} \cite{reinelt1991tsplib}, a wide set of TSP benchmark instances.
Note that most of the {\tt TSPLIB} instances have a small integrality gap for the SEP. 
The largest integrality gap is equal to 1.095, achieved by the {\tt ts225} instance, which was specifically designed to foil TSP software.
Except for that one instance, \rev{Johnson and McGeoch} \cite{johnson2007experimental} are not aware of any testbed instance with a percentage integrality gap exceeding $1.03$.

While studying theoretically the performance of the Held-Karp lower bound for the TSP, which is equivalent to the optimal lower bound of the SEP, Wolsey proved in \cite{wolsey1980heuristic} that the integrality gap for SEP is at most $\frac{3}{2}$.
Later, in his master thesis, Williamson conjectured that for metric TSP the integrality gap of subtour elimination relaxation is equal to $\frac{4}{3}$ \cite{williamson1990analysis}.
So far, this conjecture is only proved for a very special class of instances \cite{boyd2011finding}.
For cubic graphs, we know that $\frac{4}{3}$ is an upper bound for the integrality gap \cite{boyd2011tsp}.
In \cite{benoit2008finding,boyd2010structure}, Boyd and Benoit have computed the exact integrality gap for every $n \leq 12$, and, thus, they have verified by exhaustive enumeration the conjecture for all Metric TSP instances with at most 12 nodes.
In addition, Benoit and Boyd  \cite{benoit2008finding} have introduced a new class of graphs having an integrality gap that asymptotically converges to $\frac{4}{3}$ as $n$ tends to infinity.
To the best of our knowledge, nobody has been able to find, for a given number of nodes, an instance with a higher integrality gap than the one they proposed.
Curiously, we have computationally verified that those instances are not challenging for Concorde.
More recently, new families of instances with an integrality gap that asymptotically converges to $\frac{4}{3}$ were introduced in \cite{hougardy2020hard} by exploiting Tetrahedron configurations, and in \cite{HOUGARDY2014495,zhong2021lower} by exploiting three paths configurations.
\rev{All those families of instances are solvable in polynomial time}:
the instances of \cite{HOUGARDY2014495} are convex-hull-and-line TSP instances, and they are solvable in polynomial time with the algorithm proposed by \rev{Deineko} et al. \cite{deineko1994convex}. 
The instances of \cite{hougardy2020hard} can be solved using the polynomial time algorithm introduced in \cite{rubinstein2001polynomial}.
The instances of \cite{zhong2021lower} can be solved in linear time, as shown by the same author.
Nevertheless, the TSP instances introduced in \cite{zhong2021lower} are remarkable:
in practice, they are extremely challenging for Concorde.

In this paper, we introduce a new challenging integer programming problem, herein called the {\it Integer Heuri\-stic-OPT (IH-OPT)} problem, whose optimal solution provides a Metric TSP instance with a large integrality gap.
\rev{Our approach is purely computational, and we do not restrict to the generation of Euclidean or Rectilinear TSP instances as in \cite{HOUGARDY2014495,hougardy2020hard,zhong2021lower}, but we still require that the generated instance are metric}.
The decision variables of our problem are the entries of the symmetric cost matrix $\bm{C}$, while the constraints are defined by the intersection of the metric cone \cite{laurent1996graphic} with an exponential number of inequalities, one for each possible permutation of $V$, that is,
one for each possible solution of the TSP.
For example, starting from an instance $C_0$ of the {\tt TSPLIB}, taking an optimal solution of the LP relaxation of the corresponding subtour elimination model, and by solving the IH-OPT problem, we can generate a new TSP instance $C^*$ having an integrality gap larger than or equal to the integrality gap of the original instance $C_0$ (e.g., see Table \ref{tab:hardTSPLIB}).
In addition, by integrating the solution of the IH-OPT problem into a sampling procedure, we can generate several TSP instances that have a large integrality gap and are challenging in terms of runtime for Concorde (e.g., see Table \ref{tab:hardSampling}).
As a by-product of our work, we introduce the {\tt Hard-TSPLIB}, a collection of 41 small TSP instances (i.e., with $n \leq 76$) which are very challenging for Concorde in terms of runtime and number of branch-and-bound nodes.

The outline of this paper is as follows.
In Section 2, we review the background material, and we fix the notation.
Section 3 formally introduces the {\it IH-OPT} problem and discusses how it is related to previous works.
Section 4 presents the sampling procedure that we use to generate TSP instances which are hard for Concorde in terms of runtime.
In Section 5, we present our extensive computational results and describe how we generated the instances that we have included in the {\tt Hard-TSPLIB}.
Finally, in Section 6, we conclude with a discussion on future works.

\section{Background material}\label{sec:2}
In this section, we review the main formulation for the symmetric TSP and we introduce the Subtour Elimination Problem polytope.
We formally define the integrality gap for the SEP, presenting the work of \cite{benoit2008finding}, which is the foundation of our work.

An instance of the symmetric TSP can be completely defined by the symmetric matrix $\bm{C}$.
Otherwise, we can define a TSP instance using a complete undirected graph $K_n = (V, E)$ along with a cost vector $\bm{c} \in \erre_+^{\vert E \vert}$, which is given by the upper triangular matrix of $\bm{C}$ (without the diagonal).
While in this paper we focus on the general {\it Metric TSP}, and we denote it only with TSP, as defined in Section 1, two special cases that are relevant to compare our work with the literature:
\begin{enumerate}
\item The \emph{Euclidean TSP}, where the input is a collection of $n$ points in $\erre^d$, whose reciprocal distances are computed using the Euclidean norm. These are the type of TSP instances used in \cite{HOUGARDY2014495,hougardy2020hard}, with $d = 2$.
\item The \emph{Rectilinear TSP}, where the input is again a collection of $n$ points in $\erre^d$, but where the Manhattan norm (or city block or $L^1$ norm) is used to compute the distances. These are the type of TSP instances used in \cite{zhong2021lower}, with $d = 3$.
\end{enumerate}
\noindent In both cases, using the collection of points given as input, and the corresponding distance function, it is possible to define a complete graph $K_n$ with the set of nodes $V=\{1,\dots, n\}$ and the set of edges $E = V \times V$, with the corresponding cost vector $\bm c \in \erre_+^{\vert E \vert}$.
In the following, we use $\delta(S)$ with $S \subset V$, to denote the set of weighted edges $e=\{v,w\}, \ v \neq w$ with either $v\in S$, $w\not \in S$ or either $w\in S$, $v\not \in S$.
We use later the two following collection of subsets of vertices: $\mathcal{S} := \{ S \mid S \subset V,  3 \leq |S| \leq n - 3 \}$
and $\mathcal{S}_{ij} := \{ S \mid S \in \mathcal{S}, \{i, j\} \in \delta(S) \}$.

Given a connected graph $G=(V,E)$ and the cost vector $\bm{c} \in \erre_+^{\vert E \vert}$, the Travelling Salesman Problem, originally proposed in \cite{dantzig1954solution}, is formulated as follows:
\begin{align}
\label{m1:obj}  \min \quad & \sum_{e \in E} c_e x_e\\
\label{m1:degree}  \mbox{s.t.} \quad & \sum_{e\in \delta(\{v\})} x_{e} = 2 \quad \forall v\in V \quad \mbox{(Degree Constraints)}\\
\label{m1:subtour}  & \sum_{e\in \delta(S)} x_{e} \geq 2 \quad \forall S \in \mathcal{S} \quad\mbox{ (Subtour Elimination Constraints)}\\
\label{m1:bound}  & 0 \leq x_e \leq 1 \quad \;\;\; \forall e\in E \quad \quad\mbox{(Bound Constraints)}\\
\label{m1:var}  & x_e \mbox{ integer } \quad \;\;\;\; \forall e\in E.
\end{align}
\noindent The decision variable $x_e$ is equal to 1 if the edge $e=\{v,w\}$ is part of an optimal tour.
The objective function \eqref{m1:obj} minimizes the overall tour length.
Constraints \eqref{m1:degree} state that each node $v$ must have two incident edges.
The Subtour Elimination Constraints \eqref{m1:subtour} force the cut set $\delta(S)$ of every proper subset $S$ of $V$ to contain at least two edges.

If we relax the integrality constraints \eqref{m1:var}, we can define the Subtour Elimination Problem (SEP), that, given a cost vector $\bm{c}$, provides a lower bound of the optimal solution.
In the next sections, we will denote as polytope of the SEP the set 
\begin{equation}
    P_{SEP} := \{ x \in \erre^{|E|}_+ \mid \eqref{m1:degree}, \eqref{m1:subtour}, \eqref{m1:bound} \}.
\end{equation}
If we optimize the objective function \eqref{m1:obj} over $P_{SEP}$, we get exactly the Subtour Elimination Problem.

Let us consider the complete graph $K_n$.
Let us denote by $TOUR(\bm{c})$ and $SUBT(\bm{c})$, respectively the optimal value of the TSP and SEP instance defined by the cost vector $\bm{c}$.
Similarly to \cite{benoit2008finding}, we denote by $\alpha_n$ the integrality gap of SEP for $K_n$, that is, the largest possible ratio between $TOUR(\bm c)$ and $SUBT(\bm c)$:
\begin{equation}\label{eq:IG1}
    \alpha_n = \max_{\substack{\bm{c}\geq 0 \text{ is metric}}}  \dfrac{TOUR(\bm{c})}{SUBT(\bm{c})}.
\end{equation}
\noindent For $n \leq 5$, we have that $\alpha_n=1$.
The exact value for $6 \leq n \leq 12$ was computed in \cite{benoit2008finding,boyd2010structure}.
Williamson's conjecture states that $\alpha_n \leq \frac{4}{3}$ for any $n$ \cite{williamson1990analysis}.
Let's call $\alpha_n(\bm{c})$, the integrality gap of a specific TSP instance on $n$ nodes with cost vector $\bm{c}$, that is
\[ \alpha_n(\bm{c}) = \dfrac{TOUR(\bm{c})}{SUBT(\bm{c})}\]
We can thus write
\[
\alpha_n = \max_{\substack{\bm{c}\geq 0 \text{ is metric}}} \alpha_n(\bm{c})
\]

The families of instances introduced in \cite{HOUGARDY2014495,hougardy2020hard,zhong2021lower} have all the properties that $\alpha_n(\bm{c}) < \frac{4}{3}$ and $\lim_{n \rightarrow \infty} \alpha_n(\bm{c}) = \frac{4}{3}$.
However, the largest integrality gap for the instances in the {\tt TSPLIB} is equal to 1.095 \cite{reinelt1991tsplib}, which is largely less than $\frac{4}{3} = 1.\overline{333}$.

Following the idea presented in \cite{benoit2008finding}, we can divide the cost vector $\bm c$ by the optimum tour value $TOUR(\bm c)$, obtaining a new cost vector $\bm{c}'$ that still satisfies the triangle inequalities, and which leaves $\alpha_n(\bm{c})$ unchanged.
Hence, we can restrict ourselves to metric cost vectors $\bm c$ such that $TOUR(\bm{c}) = 1$, and we can transform the maximization problem into a minimization:
\begin{equation}\label{quadratic}
  \dfrac{1}{\alpha_n} = \min_{\substack{\bm{c}\geq0\text{ is metric,}\\ TOUR(\bm{c}) = 1}} SUBT(\bm{c}).
\end{equation}
Problem \eqref{quadratic} can be formulated as a mixed integer quadratic problem, where the decision variables are both the cost vector $\bm c$ and the incidence vector $\bm x$ of vertices of $P_{SEP}$.
Unfortunately, despite the recent improvements in the implementation of commercial optimization solvers, the quadratic model is intractable even for small values of $n$.
In \cite{benoit2008finding}, the authors propose a clever idea for bypassing the quadratic model using the vertex representation of $P_{SEP}$.
That is, they represent $P_{SEP}$ as the convex combinations of its vertices $\{\bm{\bar x}^{(1)},\ldots, \bm{\bar x}^{(t)}\}$, which are finitely many.
Then, for each vertex $\bm{\bar x}^{(h)}$ of $P_{SEP}$, they define an LP problem, called $\mbox{OPT}(\bar{\bm x}^{(h)})$, having an exponential number of variables and constraints.
Theoretical results in \cite{benoit2008finding,boyd2010structure} guarantee that only a subset of vertices is necessary to compute the integrality gap.
By solving the $\mbox{OPT}(\bar{\bm x}^{(h)})$ subproblem on each vertex of the previously mentioned subset, they were able to compute $\alpha_n$ for $n\leq 10$ in \cite{benoit2008finding}, and for $n \leq 12$ in \cite{boyd2010structure}.

Given the complete graph on $n$ nodes $K_n = (V,E)$, we introduce one vector  $\bm{\bar z} \in \mathbb{R}^{|E|}$ for each \rev{$\pi$} permutation of nodes, such that

\[
\bar{z}_{ij} = 
\begin{cases}
1 & \mbox{if } \pi(i) = j \mbox{ or } \pi(j) = i, \\
0 & \mbox{otherwise}.
\end{cases}
\]

Let $\mathcal{T}_n$ be the collection of the incidence vectors $\bm{\bar z} \in \erre^{|E|}$ of all the possible tours of $K_n$.
Given a vertex $\bm{\bar x}^{(h)}$ of $P_{SEP}$, the OPT$(\bar{\bm x}^{(h)})$ problem is defined in \cite{benoit2008finding} as 

\begin{align}
\label{m2:opti}   \mbox{OPT}(\bar{\bm x}^{(h)}) = \min \quad  & \sum_{\{i,j\} \in E} \bar{x}_{ij}^{(h)} c_{ij}   \\
 \mbox{s.t.} \quad
\label{m2:ctmag1}   & \sum_{\{i,j\} \in E} \bar z_{ij} c_{ij} \geq 1 & \forall \bar{\bm{z}} \in \mathcal{T}_n  \\
\label{m2:triineq}  & c_{ij} \leq c_{ik} + c_{jk} &  \forall i,j,k \in V \\
\label{m2:cnoneg}   & c_{ij} \geq 0 &  \forall \{i,j\} \in E \\
\label{m2:dual1}    & y_i + y_j - u_{ij} + \sum_{S \in \mathcal{S}_{ij}} d_S \leq c_{ij} & \forall \{i,j\} \in E \\
\label{m2:dual2}    & u_{ij} \geq 0 & \forall \{i,j\} \in E \\
\label{m2:dual3}    & d_S \geq 0 &  \forall S\in \mathcal{S} \\
\label{m2:dual1eq}    & y_i + y_j - u_{ij} + \sum_{S \in \mathcal{S}_{ij}} d_S = c_{ij} & \forall \{i,j\} \in E \mbox{ such that } \bar x_{ij}^{(h)} > 0\\
\label{m2:dual2eq}    & u_{ij} = 0 & \forall \{i,j\} \in E  \mbox{ such that } \bar{x}_{ij}^{(h)} < 1\\
\label{m2:dual3eq}    & d_S = 0 &  \forall S\in \mathcal{S}  \mbox{ such that } \sum_{ij \in \delta(S)} \bar x_{ij}^{(h)} > 2.
\end{align}

\noindent Constraints \eqref{m2:ctmag1} ensure that the optimal solution $\bm c^*$ of $\mbox{OPT}(\bar{x}^{(h)})$ for every $\bm{\bar x}^{(h)}$ is such that $TSP(\bm c^*) = 1$, as discussed in \cite{benoit2008finding}.
Herein, we call the inequalities \eqref{m2:ctmag1} the {\it TSP constraints}.
Constraints \eqref{m2:triineq} and \eqref{m2:cnoneg} ensure that the cost vectors represent a semi-metric.
Constraints \eqref{m2:dual1}--\eqref{m2:dual3} are the dual constraints associated to the dual problem of \eqref{m1:obj}--\eqref{m1:bound}.
Constraints \eqref{m2:dual1eq}--\eqref{m2:dual3eq} ensure that the vertex $\bm{\bar x}^{(h)}$ remains the optimal solution of the SEP.
If we denote by $\mathcal{C}^*$ the set of the optimal solutions of $\mbox{OPT}(\bar{\bm x}^{(h)})$, the dual slackness constraints are introduced to guarantee
\begin{equation}\label{eq:stayhome}
 \arg\min \mbox{SEP}(\bm c^*) =  \bm{x}^{(h)}, \quad \forall \ \bm c^* \in \mathcal{C}^*.
\end{equation}
Note that, in \cite{benoit2008finding,boyd2010structure}, the authors observed that, for $6 \leq n \leq 12$, $\vert \mathcal{C}^*\vert = 1$.

Hence, given the list of the vertices $\{\bm{\bar x}^{(1)},\ldots, \bm{\bar x}^{(t)}\}$ of $P_{SEP}$, the integrality gap $\alpha_n$ of $K_n$ is computed by solving
\begin{equation}\label{eq:IG2}
  \frac{1}{\alpha_n} = \min_{h = 1, \ldots, t} \mbox{OPT}(\bar{\bm x}^{(h)}) \quad\quad \Rightarrow \quad\quad \alpha_n = \max_{h = 1, \ldots, t} \frac{1}{\mbox{OPT}(\bar{\bm x}^{(h)})}.
\end{equation}
Part of the original contribution presented in \cite{benoit2008finding} is to prove that we need to consider only a subset of the $t$ vertices of $P_{SEP}$, namely, such vertices which support graphs satisfy certain properties.
We refer the reader to \cite{benoit2008finding} for the details.

In the following section, we modify the single problem $\mbox{OPT}(\bar{\bm x}^{(h)})$ in order to introduce a new NP-hard problem that we use to generate metric TSP instances with a large integrality gap.

\section{The Integer Heuristic-OPT Problem}\label{sec:IHOPT}
The goal of our work is to devise an efficient computational procedure to generate instances with a large integrality gap for the Metric TSP.
We are not interested in computing the exact value $\alpha_n$ for a fixed $n$ as in \cite{benoit2008finding}, but we focus on finding a heuristic solution to problem \eqref{eq:IG2}.

Notice that problem \eqref{m2:opti}--\eqref{m2:dual3} has an exponential number of variables due to the dual variables $d_S$.
Furthermore, the number of constraints \eqref{m2:ctmag1} is equal to the number of tours, that is $\frac{(n-1)!}{2}$, and the number of triangular inequalities constraints \eqref{m2:triineq} is $O(n^3)$.
In practice, the exact solution of $\mbox{OPT}(\bar{\bm x}^{(h)})$ is intractable even for small values of $n$.
Note that the computation of $\alpha_n$ for $n=12$ required around 24 days \cite{boyd2010structure}.
In order to solve problem \eqref{eq:IG2} heuristically, we introduce a new problem, called {\it Heuristic-OPT (H-OPT)}, which is related to problem $\mbox{OPT}(\bar{\bm x}^{(h)})$, but can be solved for larger values of $n$.
The two key ideas for introducing the new problem are:
\begin{enumerate}[(i)]
    \item To generate a TSP instance with a sufficiently large integrality gap it is unnecessary to enumerate all vertices $\{\bm{\bar x}^{(1)},\ldots, \bm{\bar x}^{(t)}\}$ of $P_{SEP}$. We can sample a subset of vertices and take the instance providing the largest integrality gap.
    \item Since we do not perform exhaustive vertex enumeration, it is unnecessary to impose the complementary slackness constraints \eqref{m2:dual1}-\eqref{m2:dual3} to force that a given vertex $\bm{\bar x}^{(h)}$ remains the same vertex of $P_{SEP}$.
    Hence, we can remove all the variables and constraints related to the slackness conditions.
\end{enumerate}
For these two reasons, given a vertex $\bm{\bar x}^{(h)} \in P_{SEP}$, we define the following LP problem:
\begin{align}
\label{m3:obj}   \mbox{H-OPT}(\bm{\bar x}^{(h)}) := \min \quad & \sum_{\{i,j\} \in E} \bar{x}_{ij}^{(h)} c_{ij}  \\
\label{m3:tour}  \mbox{s.t.} \quad &  \sum_{\{i,j\} \in E} \bar z_{ij} c_{ij} \geq 1 & \forall \bar{\bm{z}} \in \mathcal{T}_n \\
\label{m3:triineq} & c_{ij} \leq c_{ik} + c_{jk} & \forall i,j,k\in V \\
\label{m3:cemag0light} & c_{ij} \geq 0 & \forall \{i,j\} \in E.
\end{align}
\noindent In practice, we have relaxed the problem $\mbox{OPT}(\bar{\bm x}^{(h)})$ by removing constraints \eqref{m2:dual1}--\eqref{m2:dual3}.
Notice that in H-OPT we have only $|E|$ cost variables, $\frac{(n-1)!}{2}$ TSP constraints \eqref{m3:tour} (one for each tour), and $O(n^3)$ triangles inequalities \eqref{m3:triineq} that define the metric cone \cite{laurent1996graphic}.
We can solve this LP problem by cutting planes, by separating both families of constraints.

\rev{A critical point in the solution of the H-OPT problem by branch-and-cut is the separation of (maximally violated) TSP constraints \eqref{m3:tour}. The \emph{separation problem} SP is defined as follows}. Given a cost vector $\bm{\bar c} \in \erre^{|E|}_+$, we look for a tour whose incidence vector $\bm{\bar z}$ verify the following:
\begin{equation}
  \sum_{\{i,j\} \in E} \bar z_{ij} \bar{c}_{ij} < 1.
\end{equation}
\noindent Notice that any tour that satisfies the previous relation gives a violated TSP constraint.
However, to prove that no tour violates the TSP constraint, we need to verify the following:
\begin{equation}\label{eq:tspsep}
  \min_{\bm z \in \mathcal{T}_n} \left\{ \sum_{\{i,j\} \in E} z_{ij} \bar{c}_{ij} \right\} \geq 1.
\end{equation}
Thus, to add a new TSP constraint, we solve a TSP instance for a specific cost vector $\bm{\bar c}$.
\noindent In our implementation, we separate TSP constraints by first solving the TSP instance given by $\bm c^*$ using the LK-H heuristic \cite{helsgaun2000effective,helsgaun2009general}, and whenever the heuristic fails to find a violated tour, we solve \eqref{eq:tspsep} by embedding Concorde in the code.
Section 5.5 describes the details of our implementation.
Clearly, the TSP constraints make this problem challenging, as stated in the following lemma.

\begin{lemma}\label{lemma1}
The H-OPT problem is NP-hard.
\end{lemma}

First of all, we prove the following Lemma.
\begin{lemma}\label{lemma:sepopt}
\rev{Let $n=|V|$. If \mbox{H-OPT} can be solved in $n^{O(1)}$ time, then the separation problem SP can be solved in $n^{O(1)}$ time.}
\end{lemma}
\begin{proof}[Sketch]
\rev{ Let $P$ be the polyhedron defined by equations \eqref{m3:tour}--\eqref{m3:cemag0light}}. The proof is implicit in the proof of a well-known theorem of Gr{\"o}tschel  et al. \cite{GroetschelLovaszSchrijver1981} which states that for well-described \rev{polyhedron} $P$ the optimization problem can be solved in polynomial time if and only if the \rev{separation problem SP} can be solved in polynomial time.
Note that the requirement of well-described polyhedral $P$ is only used for obtaining a time complexity polynomial in the input size of the problem. For example, if the length $L$ of the input needed to describe the polyhedral $P$ does not satisfy $n\leq L$, then even if the separation could be resolved in time $n^{O(1)}$ we could not claim that the \rev{separation problem SP} can be resolved in polynomial time.

To see how the proof of this lemma is implicit in \cite{GroetschelLovaszSchrijver1981}, it is sufficient to note that  in \cite{GroetschelLovaszSchrijver1981} it is not asked to provide the input size as tight as possible. 
Actually, it is possible to make the encoding dimension $L$ of a polyhedron $P$ greater than $n$ by adding, for example, a string with $n$ zeros.
This ensures $n\leq L$, and the claim follows.
\end{proof}

\begin{proof}[Lemma \ref{lemma1}]
The claim follows by showing that the Hamiltonian Cycle problem (HC) can be solved in polynomial time by an oracle machine with an oracle for H-OPT. More precisely, we reduce in polynomial time HC to the separation problem SP of H-OPT. Then, by Lemma \ref{lemma:sepopt} an oracle machine for H-OPT that takes $n^{O(1)}$ time implies that the separation problem SP (and therefore HC) can be solved in $n^{O(1)}$ time. 

Consider the following sets:

\begin{eqnarray}
    M &=& \{ \bm{c} \in \mathbf{R}^m \; \vert \; c_e > 0,  \forall e\in [m] \, c_{ij} \leq c_{ik} + c_{jk} \; \forall i,j,k \in [n]\}, \\
    TSP &=& \{ \bm{c} \in \mathbf{R}^m \; \vert, \bm{c}^T\bm{x} \geq 1 \; \forall \bm{x} \in \mathcal{T} \}, \\
    P &=& M \cap TSP.
\end{eqnarray}

Let $G=(V,E)$ be an undirected graph on $n$ nodes, that is $\vert V \vert = n$.
This graph defines an instance of the HC problem, namely the decision problem that searches for a cycle in a graph.
We can perform a standard polynomial time reduction to get a Metric TSP and then normalize costs as follows:

\begin{equation}
    \hat{c}_e = \begin{cases} (1 - \frac{\epsilon}{2})/n & e \in E, \\ (2-\epsilon)/n & \mbox{otherwise,} \end{cases}
    \label{eq:poly-transf}
\end{equation}
where $1 > \beta(n) > \epsilon > 0$, and the definition of $\beta(n)$ will be clarified later in the proof of this theorem.

Note that this instance of the TSP is metric: identity and symmetry are obvious, and the triangle inequalities can be verified case-by-case.
Note also that, the graph contains a Hamiltonian cycle if and only of the optimal solution of the TSP is $1-\frac{\epsilon}{2} < 1$, and thus such cost vector is in $M \setminus TSP$.
On the opposite side, if the graph does not admit Hamiltonian cycle, then the TSP solution must contain at least one edge of length $(2-\epsilon)/n$.
Then, the value of the optimal tour would be at least
\begin{equation}
    \dfrac{2 - \epsilon}{n} + \dfrac{(n-1)(1-\frac{\epsilon}{2})}{n}
    =  \frac{1}{n} + 1 - \frac{\epsilon}{n}\left(\frac{1}{2} + \frac{n}{2} \right).
    \label{eq:hc_value}
\end{equation}

We observe that if $\varepsilon < \dfrac{2}{1+n}$, then \eqref{eq:hc_value} is greater than 1. Thus, we can set $\beta(n) =  \dfrac{2}{(n+1)}$.
Note that from now on we have proved that a graph $G = (V, E)$ on $n$ nodes admits a Hamiltonian cycle if and only if the solution of the associated TSP with cost vector $\hat{\bm{c}}$ is less than 1. 
By putting everything together, we see that the Hamiltonian Cycle problem (HC) can be solved in polynomial time by an oracle machine with an oracle for H-OPT as follows:
\begin{enumerate}
    \item Let $G= (V, E)$ be an undirected graph on $n$ nodes.
    \item Use equation \eqref{eq:poly-transf} to obtain a TSP instance with cost vector $\hat{\bm{c}}$.
    \item By Lemma \ref{lemma:sepopt}, an oracle machine for H-OPT that takes $n^{O(1)}$ time implies that we can decide \rev{whether} $\hat{\bm{c}}\in M$ in $n^{O(1)}$ time. If $\hat{\bm{c}}\in M$ then $G$ does not admit \rev{a Hamiltonian} cycle. Otherwise, it admits \rev{a Hamiltonian} cycle.
\end{enumerate}
\qed
\end{proof}

\rev{We remark that our proof only yields to NP-hardness under Turing reductions, as we have shown that there is an NP-complete decision problem, namely HC, that can be Turing-reduced to H-OPT.
However, this also implies ``hardness'' for our problem, as an existence of a polynomial time algorithm for H-OPT would imply $P=NP$.}

Note also that since we have removed from our problem the dual slackness constraints, the relation \eqref{eq:stayhome} is not necessarily satisfied when $\bm c^*$ is the optimal solution of H-OPT$(\bar{\bm x}^{(h)})$.
However, we can prove the following lemma which states that the solution of the H-OPT$(\bar{\bm x}^{(h)})$ problem provides a cost vector $\bm c^*$ corresponding to a TSP instance with an integrality gap greater than or equal to any instance yielding $\bar{\bm x}^{(h)}$.

\begin{lemma}\label{lem:gap}
Let us consider a TSP instance $\bm{c}_0$ such that $TOUR(\bm{c}_0) = 1$, and let $\bm{\bar{x}}^{(0)}$ be an optimal solution of $SEP(\bm{c}_0)$.
We define a second TSP instance by the cost vector
\[
  \bm{c}_1 = \arg\min \mbox{H-OPT}(\bm{\bar{x}}^{(0)}).
\]
Then, the following relation holds
\begin{equation}
  \frac{TOUR(\bm{c}_1)}{SUBT(\bm{c}_1)} \geq \frac{TOUR(\bm{c}_0)}{SUBT(\bm{c}_0)}.
\end{equation}
\end{lemma}
\begin{proof}
Since $\bm{\bar{x}}^{(0)}$ is a feasible solution of the SEP, we have that
\[
  SUBT(\bm{c}_1) \leq \bm{c}_1^T \bm{\bar{x}}^{(0)}
\]
By definition of H-OPT$(\bm{\bar{x}}^{(h)})$, we have $\bm{c}_1$ is the cost vector that realizes the minimum of $\bm{c}^T \bm{\bar{x}}^{(h)}$ among all the metric vectors such that $TOUR(\bm{c})=1$.
By hypothesis, $TOUR(\bm{c}_0) = 1$ and this implies, $\bm{c}_0 \bm{z} \geq 1$ for all $\bm{z}$ 0-1 incidence vector.
Thus, $\bm{c}_0$ is feasible for H-OPT$(\bm{\bar{x}}^{(0)})$ and it holds
\[
\bm{c}_1^T \bm{\bar{x}}^{(0)} \leq \bm{c}_0^T \bm{\bar{x}}^{(0)} = SUBT(\bm{c}_0),
\]
where the last equation holds by definition.
Thus,
\[
  \dfrac{TOUR(\bm{c}_0)}{SUBT(\bm{c}_0)} = \dfrac{1}{SUBT(\bm{c}_0)}
  \leq
  \dfrac{1}{\bm{c}_1^T \bm{\bar{x}}^{(0)}}
  \leq
  \dfrac{1}{SUBT(\bm{c}_1)} = \dfrac{TOUR(\bm{c}_1)}{SUBT(\bm{c}_1)}.
\]
\qed
\end{proof}

\paragraph{From large integrality gaps to hard instances.}
The TSP solver Concorde \cite{applegate2006traveling} and the LK-H heuristic that we use for separating TSP constraints \cite{helsgaun2000effective} can handle only integer costs.
However, the H-OPT problem generates fractional cost vectors.
Hence, we need to devise a method to transform the fractional costs into integer values, by, for example, multiplying for a large constant $\tau$ and then rounding to the nearest integer, that is, $c_{ij} = \mbox{round}(\tau c^*_{ij})$.
Using this cost transformation, we cannot guarantee that the integrality gap for $c_{ij}$ remains the same of $c^*_{ij}$.

The TSP constraints \eqref{m3:tour} have a right-hand side equal to 1 because we have divided the cost vector $\bm c$ by the minimum tour length $TOUR(\bm c)$.
If we divide the cost vector by the quantity $\left( \frac{TOUR(\bm c)}{\Delta} \right)$, where $\Delta>0$ is a large positive constant,
we have to change the right hand side of \eqref{m3:tour} to $\Delta$, while still getting an equivalent LP problem.
If $\Delta$ is large enough (see Section 5.5), we can also add the integrality constraint on the variable $c_{ij}$. 
In practice, in order to generate integer cost vectors, we have introduced the following ILP problem:
\begin{align}
  \mbox{IH-OPT}(\bm{x}^{(h)}) := \min \quad & \sum_{\{i,j\} \in E} \bar{x}_{ij}^{(h)} c_{ij}  \\
  \mbox{s.t.} \quad & \eqref{m3:triineq}, \eqref{m3:cemag0light} \\
\label{m4:delta}  & \sum_{\{i,j\} \in E} \bar z_{ij} c_{ij} \geq \Delta & \forall \bar{\bm{z}} \in \mathcal{T}_n \\
  & c_{ij} \mbox{ integer } & \forall \{i,j\} \in E.
\end{align}

Clearly, IH-OPT is very challenging, as it is an integer program with as many constraints as TSP tours plus the number of triangle inequalities. 
We \rev{have computational evidence} that the integrality gap of the integer problem could be smaller than those obtained by solving H-OPT.
Surprisingly, our computational results show that the TSP instances obtained while solving the IH-OPT problem are very challenging for Concorde.

In the next section, we show how we use the IH-OPT problem to search for very challenging instances for Concorde.

\section{A sampling procedure for generating hard instances}
\label{sec:sampling}

A standard procedure for searching heuristically for (suboptimal) solutions of an optimization problem is based on uniformly sampling points (i.e., solutions) of the feasible region \cite{bertsimas2004solving}.
In our context, to find a heuristic solution for problem \eqref{eq:IG2}, we could sample a fixed number of vertices of $P_{SEP}$ and retain the vertex yielding the minimum value for IH-OPT, that is, yielding the largest integrality gap. 
However, directly sampling the vertices of the $P_{SEP}$ is impractical, due to its exponential number of subtour constraints.
We could instead easily sample random cost vectors and generate vertices of $P_{SEP}$ by solving directly problem \eqref{m1:obj}--\eqref{m1:bound}.

A possibility for generating random cost vectors consists of generating $n$ random points in a Euclidean space, and then computing all the pairwise distances using a given distance (e.g., a distance induced by the Minkowski norm).
As observed from our preliminary tests, this procedure leads very often to an integral, and hence useless, vertex of $P_{SEP}$. 
In practice, this procedure takes a long time before returning a fractional vertex of $P_{SEP}$.
Furthermore, it only samples instances of the Euclidean TSP, implicitly excluding some remarkable metric TSP instances, such as the ones provided in \cite{benoit2008finding}.

To generate a random cost vector, we have designed a different approach.
First, we sample a random point within the metric polytope \cite{laurent1996graphic} using the {\it hit-and-run} algorithm \cite{smith1984efficient}, which generates a random metric TSP instances $\bm{c} \in \erre_+^ {|E|}$.
Second, we get a (random) vertex $\bm{\bar x}^{(h)}$ by solving SEP$(\bm c)$ via the simplex algorithm.
Since the cost vector is a uniformly random point of the metric polytope, we expect that $\bm{\bar x}^{(h)}$ is a random vertex of $P_{SEP}$. 
We also expect a good variety among the sampled vertices: for instance, with $n = 15$, after the sampling of 999 vertices, it took only less than three seconds to find one vertex both non-integer and not yet sampled.
In the next paragraphs, we briefly review the hit-and-run algorithm, and we detail how we use the sampling procedure to generate hard metric TSP instances.

\paragraph{The hit-and-run algorithm}
\rev{The hit-and-run algorithm \cite{smith1984efficient} is designed to sample from a bounded set $P$ uniformly.}
The basic steps of the algorithm are:
\begin{enumerate}
\item Pick a point $\bm{x}_k \in P \subset \erre^m$.
\item Generate a random direction $\bm{d}_k$ uniformly distributed over the unitary hyper-sphere centered in  $\bm{x}_k$.
\item Generate a random point $\bm{x}_{k+1} = \bm{x}_{k} + \lambda \bm{d}_k$ uniformly distributed over the \emph{line set} $\{\bm{x} \in P \ \vert \ \bm{x} = \bm{x}_k + \lambda \bm{d}_k, \ \lambda \in \erre  \}$.
\item If a stopping criterion is met, stop (e.g., terminate after a fixed number of iterations). Otherwise, $\bm{x}_{k} \gets \bm{x}_{k+1}$ and repeat from Step 2.
\end{enumerate}

The points sampled with this procedure converge in total variation to a uniform distribution, as proven in \cite{smith1984efficient}.
The time required to have a sample that effectively \rev{approximates} the uniform distribution is polynomial in the dimension, as proven in \cite{lovasz1999hit}.
Later, a modification of this algorithm to sample \rev{a set of points that converges} to an arbitrary target distribution was introduced in \cite{romeijn1994simulated}.
The interested reader can find an in-depth review of hit-and-run algorithms in \cite{zabinsky2013hit}.

\paragraph{The random sampling algorithm}
As noted in Section \ref{sec:2}, the constraints \eqref{m2:triineq}--\eqref{m2:cnoneg} define the metric cone:
\begin{equation}
  C_{MET} := \{\bm c \in \erre^{|E|}_+  \mid c_{ij} - c_{ik} - c_{jk} \leq 0, \forall i,j,k \in V \}.
\end{equation}
\noindent If we add the perimeter inequality (also known as the {\it homogeneous} triangle inequality) we get the metric polytope
\begin{equation}
  P_{MET} := C_{MET} \cap \{ c_{ij} + c_{ik} + c_{jk} \leq 2, \forall i,j,k \in V  \}.
\end{equation}
\noindent
Note that for every metric $\bm{c} \in C_{MET}$, there exist an $\gamma$ such that $\gamma \bm{c} \in P_{MET}$ \cite{laurent1996graphic}.
For this reason, without loss of generality, we can sample from $P_{MET}$ instead of $C_{MET}$.
\rev{Sampling from $P_{MET}$ instead of $C_{MET}$ guarantees to remain in the framework described in \cite{lovasz1999hit}, since $P_{MET}$ is a convex compact set and, hence, a convex body.}

Algorithm \ref{alg:herSEP} presents our procedure for sampling metric TSP instances from the set ${P}_{MET}$.
Given the size of the metric space $m=|E|$ and the number of required sampled points $r$, the algorithm initialized in Step 1 the empty set $R$.
Then, until the number of sampled points is equal to $r$, steps 3--7 are iterated.
Step 3 generates a uniformly distributed random cost vector $\bm c$ from the open set ${P}_{MET}$ using the hit-and-run algorithm.
Step 4 generates an optimal vertex $\bm x$ for SEP$(\bm{c})$, $\bm{x}\in P_{SEP}$.
If the vertex $\bm x$ is fractional and does not belong to $R$, it is added to $R$.
Finally, the procedure returns the set $R$ of $r$ vertices of $P_{SEP}$ of dimension $m$.

\begin{algorithm}[t!]
\DontPrintSemicolon
\SetAlgoLined
 \KwInput{$m = |E|$, the size of the TSP instance}
 \KwInput{$r$, the required number of vertices}
 \KwOutput{$R$, a collection of vertices of $P_{SEP}$}

 $R \gets \emptyset$\;
 \While{$\vert R \vert < r$}{
    $\bm{c} \gets HitAndRun({P}_{MET})$\;
    $\bm{x} \gets \arg\min \{ SEP(\bm{c}) \}$\;
    \If{$\bm{x}$ is fractional $\mbox{ {\bf and} } \bm{x} \notin R$}{
      $R \gets R \cup \{\bm x\}$\;
    }
 }
 \Return R
\caption{Sampling vertices of the $P_{SEP}$ by using the metric polytope.\label{alg:herSEP}}
\end{algorithm}

\section{Generating the {\tt Hard-TSPLIB}}\label{HardTSPLIB}

The objective of our computational experiments is to generate the {\tt Hard-TSPLIB}, a collection of small metric TSP instances having a large integrality gap and are challenging for the Concorde solver.
The instances of the {\tt Hard-TSPLIB} are generated by using a branch-and-cut solver for problem $\mbox{IH-OPT}(\bm{\bar{x}}^{(h)})$.

In the following paragraphs, first, we present the hard instances generated starting from the {\tt TSPLIB}.
Second, we present the hard instances generated using the random sampling procedure presented in Section \ref{sec:sampling}.
 Third, we compare the runtime of Concorde for solving our hard instances with the runtime required for solving the \emph{Rectilinear 3-Dimensional instances} \cite{zhong2021lower}.
In this work, the author compares the introduced instances with other works, namely \cite{tsplib,benoit2008finding,HOUGARDY2014495,hougardy2020hard}, showing that author's instances are computatinoally harder.
Thus, we compare ourselves with \cite{zhong2021lower} and deduce other comparisons from the context.
Then, we visually analyze the structure of the hard instances we have generated.
Finally, we discuss the implementation details of our algorithms.

Note that all the computations described in the following paragraphs are executed on a single node of an HPC cluster running CentOS, having an Intel CPU with 32 physical cores working at 2.1 GHz, and 64 GB of RAM.
We compiled our code with the GNU C++ compiler v8.3, with the following flags {\tt -O2 -D\_REENTRANT -m64 -ffast-math -DNDEBUG -Wall -march=native}.

\subsection{Implementation details}
Our computational procedure for generating hard metric TSP instances has two core algorithms:
\begin{enumerate}
  \item The branch-and-cut algorithm for solving problem IH-OPT$(\bm{\bar{x}}^{(h)})$ (see Section 3).
  \item The sampling procedure from the metric polytope (see Section 4).
\end{enumerate}
In the following paragraphs, we first describe our implementation of the two algorithms.

\paragraph{Solving IH-OPT by branch-and-cut.}
The problem IH-OPT$(\bm{\bar{x}}^{(h)})$ is solved by branch-and-cut by using the Gurobi commercial solver and using the C++ programming language.
In our implementation, we dynamically add both the triangle inequalities and the TSP constraints, in order to keep the core LP problem as small as possible.
Note that a complete enumeration of triangle inequalities could exhaust the memory of a standard computer for medium values of $n$.

For the separation of triangular inequalities, we follow the strategies used in \cite{grotschel1989cutting}:
%
at each iteration, we add a fixed number $k$ of the most violated triangular inequalities, until no more violated triangle inequalities exist.
The separation of triangle inequalities is carried over both during the solution of the LP relaxation and when encountering a new incumbent integer solution.

The separation of TSP constraints is more challenging since it corresponds to the solution of a TSP instance, as shown in \eqref{eq:tspsep}.
As we only need a violated cut, we first separate the TSP constraints heuristically by solving the TSP instance using the Lin-Kernigan local search procedure \cite{lin1973effective}, as implemented in Concorde\footnote{\url{http://www.math.uwaterloo.ca/tsp/concorde/downloads/codes/src/co031219.tgz}.}.
Note that the LK-H heuristic only handles integer costs.
Hence, we use the support of Gurobi for lazy constraints, which allows running our separation procedure only on incumbent integer solutions.
Whenever the LK-H heuristic returns a TSP solution whose incidence vector $\bm z$ satisfies $\sum_{\{i,j\}\in E }z_{ij} c_{ij} \geq  \Delta$, with $\Delta$ as described in Section \ref{sec:IHOPT}, we run a second exact TSP algorithm.
For the exact separation of TSP constraints we use Concorde compiled using CPLEX 12.8 as LP solver, by using the default parameter settings.

The optimal solution of the LK-H heuristic is given as warm start to the TSP branch-and-cut implementation.

Notice that before starting the solution of the integer problem IH-OPT$(\bm{\bar{x}}^{(h)})$, we solve its LP relaxation H-OPT$(\bm{\bar{x}}^{(h)})$ via cutting planes by only separating the TSP constraints with the LK-H heuristic (see, for instance, the computational results in Table \ref{tab:2.1}).
Once the LK-H heuristic does not find any TSP violated constraint, we stop, and we collect all the triangle inequalities violated by no more than a threshold $\tau=0.05$, and all the generated TSP constraints, to initialize the first pool of cut for our branch-and-cut algorithm.

\paragraph{Sampling the metric polytope by hit-and-run.}
The sampling procedure described in Algorithm \ref{alg:herSEP} is implemented in {\tt Python 3.8.2}.
For the hit-and-run algorithm at Step 3, we have used the implementation provided by \cite{her}, which is based on the original algorithm introduced in \cite{smith1984efficient}.
The solution of the SEP problem in Step 4 is implemented using the python wrapper of {\tt Gurobi}.
Since the sampling procedure is very fast, we did not port this procedure to C++.

\paragraph{Parameters tuning.}
In our implementation, we had only to decide a value for the parameter $\Delta$ in \eqref{m4:delta}.
If we only consider the H-OPT formulation, namely the formulation without the integer costs, we can imagine to multiply each cost for a fixed quantity $\omega$ and obtain a tour that costs $\omega\Delta$.
However, once we move to integer costs, the parameter $\Delta$ becomes more important because it is strongly related to the number of values that can be taken as cost coefficients.
For instance, if $\Delta = n$, \rev{the solution having all $c_{ij}$ equal to 1 is optimal.}
This does not lead to an hard-to-solve instance.
On the contrary, if $\Delta$ is too big, we miss such ``degeneracy'' of the costs on some edges.
In practice, we have observed that for different values of $\Delta$ we get TSP instances of different (runtime) difficulty.
Table \ref{tab:ccrossval} shows the computational results for solving IH-OPT over sampled TSP instances using different values of $\Delta$.
Since the hardest instances were generated while using $\Delta=1000$, we fix this value in all of the tests reported in this paper.
In future work, we plan to investigate further the impact of the parameter $\Delta$.

\begin{table}[t!]
\centering
\caption{Impact of the parameter $\Delta$ on the instances generated by IH-OPT. }
\label{tab:ccrossval}       
\begin{tabular}{lrrr}
\hline\noalign{\smallskip}
& $\Delta$ = 100 & $\Delta$ = 1000 &  $\Delta$ = 10000\\
\noalign{\smallskip}\hline\noalign{\smallskip}
\textbf{gr24, hard} & & & \\
{Integrality gap} & 1.146 & 1.220 & 1.220\\
{avg. Concorde mean time}  & 3.263 & 16.3 & 10.77\\
{std. dev Concorde time} & 0.977 & 2.0 & 1.924\\
\textbf{bayg29, hard} & & & \\
{Integrality gap} & 1.159 & 1.186 & 1.187\\
{avg. Concorde mean time}  & 4.157 & 82.5 & 32.85\\
{std. dev Concorde time} & 1.271 &  17.6 & 9.200\\
\textbf{bays29, hard} & & & \\
{Integrality gap} & 1.197 & 1.229 & 1.228\\
{avg. Concorde mean time}  & 11.108 & 64.2 & 47.392\\
{std. dev Concorde time} & 1.790 & 18.2 & 16.137\\
\noalign{\smallskip}\hline\noalign{\smallskip}
\end{tabular}
\end{table}

\subsection{Generating hard instances from the {\tt TSPLIB}}
The {\tt TSPLIB} contains 20 instances with less than 76 nodes.
Only 13 of them have a fractional solution for the SEP (i.e., the integrality gap is greater than 1).
For 12 of these 13 instances, we have generated a corresponding hard instance by solving H-OPT$(\bm{\bar{x}}^{(h)})$, where $\bm{\bar{x}}^{(h)}$ is the fractional solution of SEP solved via the simplex algorithm.

If we denote by $\bm c_0$ the cost vector of the {\tt TSPLIB} instance, and by $\bm c^*$ the optimal solution of H-OPT$(\bm{\bar{x}}^{(h)})$, by Lemma \ref{lem:gap}, we have
\begin{equation*}
  \frac{TOUR(\bm{c}^*)}{SUBT(\bm{c}^*)} \geq \frac{TOUR(\bm{c}_0)}{SUBT(\bm{c}_0)},
\end{equation*}
that is, the TSP instance $\bm c^*$ has an integrality gap larger than or equal to $\bm c_0$.
If instead, we solve the problem IH-OPT$(\bm{\bar{x}}^{(h)})$, we cannot guarantee the previous relation, but in practice, we get more challenging instances with almost the same integrality gap.
For this reason, all the following results are obtained by solving the IH-OPT problem.
Table \ref{tab:hardTSPLIB} reports the detailed results for the generation of hard instances from the {\tt TSPLIB}.
The table first reports the name and the dimension $n=|V|$ of the original instance.
The third and fourth columns report the integrality gap of $\bm c_0$ (easy instance) and $\bm c^*$ (hard instance).
In the remaining six columns, the table shows the average runtime in seconds (with the standard deviations), and the average number of branch and bound nodes for solving with Concorde first the {\tt TSPLIB} instance, and later the corresponding {\tt Hard-TSPLIB} instance.
The averages are computed over 10 independent runs of Concorde, using 10 different seeds.
Finally, the last column reports whether the optimal SEP solution of $\bm c_0$ is the same optimal solution of $\bm c^*$.

The results of Table \ref{tab:hardTSPLIB} show that the small {\tt TSPLIB} instances have a very small integrality gap and are extremely easy for Concorde.
They are solved within a fraction of seconds at the root node of the branch-and-cut tree.
On the contrary, the {\tt Hard-TSPLIB} have a significantly larger integrality gap, and they require, on average, several seconds (or up to several hours) to be solved to optimality by Concorde.
As expected, a larger integrality gap at the root node implies a larger branch-and-cut tree, that is, a larger number of nodes (column `BC n.').
However, this is not always true, because among the three instances with 48 nodes ({\tt att48}, {\tt gr48}, and {\tt hk48}), the instance with the smaller integrality gap requires the largest number of branch-and-cut nodes.
Indeed, visiting a larger search tree implies a longer runtime.
We remark that the instance {\tt brazil58\_hard} is not solved by Concorde within a timeout of 24 hours.

\begin{table}[t!]
\caption{{\tt Hard-TSPLIB} instances generated from the {\tt TSPLIB}. The columns give the instance name, number of nodes $|V|$, the integrality gap for the TSPLIB instance $\bm c_0$ and the HardTSPLIB $\bm c^*$; average runtime (and standard deviation) and number of BC nodes over 5 independent runs of Concorde, compiled with CPLEX 12.8. The last column reports if equation \eqref{eq:stayhome} holds. The last three instances marked with (*) are solved only once due to the large running time. The instance {\tt brazil58} reached a timeout of 24 hours.}
\label{tab:hardTSPLIB}       
\begin{adjustbox}{width=\textwidth}
\begin{tabular}{lcrrrrrrrrc}

\hline\noalign{\smallskip}
& & \multicolumn{2}{c}{{ $\alpha_n(\bm{c})$}} & \multicolumn{3}{c}{{TSPLIB} - $\bm c_0$} & \multicolumn{3}{c}{{\sc HardTSPLIB} - $\bm c^*$} & Does Eq \eqref{eq:stayhome}\\
name & $|V|$   & $\bm{c}_0$ & {$\bm{c}^\ast$ } & runtime & stddev & BC n.  & runtime & stddev & BC n. & hold?\\
\noalign{\smallskip}\hline\noalign{\smallskip}
gr24          & 24      & 1.000   & 1.220       & 0.028        & 0.004    & 1       & 16.3   & 2.0     & 33.0        & \xmark    \\
bayg29        & 29      & 1.001   & 1.186       & 0.030        & 0.004    & 1       & 82.5   & 17.6    & 158.6        & \xmark    \\
bays29        & 29      & 1.003   & 1.229       & 0.030        & 0.004    & 1       & 64.2   & 18.2    & 106.2        & \xmark   \\
dantzing42    & 42      & 1.003   & 1.150       & 0.031        & 0.005    & 1       & 66.4   & 51.7    & 123.0        & \cmark    \\
swiss42       & 42      & 1.001   & 1.140       & 0.031        & 0.005    & 1       & 1466.3 & 1122.5  & 369.8        & \xmark   \\
att48         & 48      & 1.002   & 1.135       & 0.031        & 0.005    & 1       & 1204.4 & 444.6   & 2876.2       & \xmark   \\
gr48          & 48      & 1.018   & 1.222       & 0.031        & 0.005    & 1       & 97.4   & 18.3    & 132.6        & \xmark   \\
hk48          & 48      & 1.001   & 1.215       & 0.031        & 0.005    & 1       & 185.0  & 55.9    & 265.8        & \xmark   \\
eil51         & 51      & 1.008   & 1.285       & 0.032        & 0.006    & 1       & 743.0  & 148.5   & 1253.8       & \xmark \\        
brazil58         & 58 & 1.002 & 1.163 &	0.033 & 0.005 & 1 & {\it 86400.0} & (*) & {\it 55139} & \xmark  \\
st70         & 70      &  1.006	& 1.280 & 0.033 &	0.005 &	1 & 37759.1 & (*) & 28227 &  \xmark \\        
pr76         & 76   &   1.029 & 1.282 & 0.035 &	0.007 &	1 & 63181.6 & (*) & 37459 & \cmark  \\        
\noalign{\smallskip}\hline
\end{tabular}
\end{adjustbox}
\end{table}

Table \ref{tab:2.1} and Table \ref{tab:2.2} report the computational results for generating the {\tt Hard-TSPLIB} instances while solving H-OPT$(\bm{\bar{x}}^{(h)})$ and IH-OPT$(\bm{\bar{x}}^{(h)})$, respectively.
For each instance, the table reports the number of cuts generated and the runtime for the triangle inequalities separation ({\it trian.}), the TSP constraints separated by the LKH heuristic ({\it LKH-cuts}), and the TSP constraints separated by Concorde ({\it TSP-cuts}).
For the solution of IH-OPT$(\bm{\bar{x}}^{(h)})$, the table gives also the total number of branch-and-bound nodes, the lower bounds (LB), and the upper bounds (UB): when the LB and UB are equal the instance is solved to optimality.
Concerning the separation algorithms, the violated triangle inequalities are identified in a very short time, and they have almost no impact on the overall runtime.
The heuristic separation of TSP constraints using the LKH heuristic is very effective, but the runtime begins to be important.
The exact separation of TSP constraints is one of the two runtime bottlenecks for the generation of hard instances. 
For example, for the instance {\tt swiss42}, most of the time is spent on the exact separation of a TSP cut. 
Finally, notice that for the instances {\tt pr76, eil76, rat99, kroB100, kroC100}, the solution by branch-and-cut of IH-OPT$(\bm{\bar{x}}^{(h)})$ hits the time limit of 24 hours (86400 seconds).
Among those instances, only for {\tt pr76}, we do generate a hard TSP instances; in all other cases, we were not able to find an integer cost vector satisfying all triangle inequalities and all TSP constraints, that is, an optimal integer solution for problem IH-OPT$(\bm{\bar{x}}^{(h)})$.

\begin{table}[]
\centering
\caption{Computational results for the heuristic solution of the LP problem H-OPT$(\bm{\bar{x}}^{(h)})$ by cutting planes. At this stage, we only separate the TSP cuts with the LK-H heuristic. Columns 3 and 4 report the number of triangular inequality and TSP constraints added. Columns 5 and 6 report the runtime for separating those inequalities. The last column reports the overall runtime in seconds.}
\label{tab:2.1}
\begin{tabular}{llrrrrr}
\noalign{\smallskip}\hline
                  &               & \multicolumn{2}{c}{{\sc  cuts}}                                 & \multicolumn{3}{c}{\sc runtime}                                                         \\
\multicolumn{1}{c}{name} & \multicolumn{1}{c}{$\vert V \vert$} & \multicolumn{1}{c}{trian.} & \multicolumn{1}{c}{LK-H} & \multicolumn{1}{c}{trian.} & \multicolumn{1}{c}{LK-H} & \multicolumn{1}{c}{total} \\ \hline
gr24                         & 24                     & 1500                       & 61                           & 0.00                       & 1.1                         & 2.2                       \\
bayg29                       & 29                     & 2394                       & 66                           & 0.00                       & 3.1                         & 6.5                       \\
bays29                       & 29                     & 2318                       & 70                           & 0.00                       & 2.5                         & 5.2                       \\
dantzig42                    & 42                     & 6515                       & 88                           & 0.00                       & 7.2                         & 38.7                      \\
swiss42                      & 42                     & 5083                       & 87                           & 0.00                       & 15.0                        & 46.5                      \\
gr48                         & 48                     & 10011                      & 182                          & 0.00                       & 20.7                        & 93.6                      \\
hk48                         & 48                     & 7494                       & 177                          & 0.00                       & 42.2                        & 78.7                      \\
eil51                        & 51                     & 7358                       & 123                          & 0.00                       & 25.8                        & 69.1                      \\
brazil58                     & 58                     & 9271                       & 136                          & 0.00                       & 85.9                        & 287.5                     \\
st70                         & 70                     & 21914                      & 411                          & 0.00                       & 116.7                       & 1593.4                    \\
eil76                         & 76                     & 27885                      & 418                          & 0.00                       & 206.3                       & 2684.5                    \\
pr76                         & 76                     & 26233                      & 473                          & 0.00                       & 141.7                       & 4716.2                    \\
rat99                        & 99                     & 61300                      & 757                          & 0.00                       & 848.0                       & 190715.4                  \\
kroB100                      & 100                    & 60300                      & 883                          & 0.00                       & 841.5                       & 162504.0                  \\
kroC100                      & 100                    & 62501                      & 1084                         & 0.01                       & 757.4                       & 131369.2 \\
\noalign{\smallskip}\hline
\end{tabular}
\end{table}

\begin{table}[]
\centering
\caption{Computational results for the solution of the integer problem IH-OPT$(\bm{\bar{x}}^{(h)})$ by branch-and-cut. Columns 3, 4, and 5 report the number of cuts added for each type. Columns 6, 7, and 8  report the runtime of the separation algorithms for triangular, LKH and exact TSP cuts. Column 9 reports the runtime and column 10 the number of branch-and-bound nodes required. The last two columns gives the lower (LB) and upper bounds (UB).}
\label{tab:2.2}
\begin{adjustbox}{width=\textwidth}
\begin{tabular}{llrrrrrrrrrr}
\noalign{\smallskip}\hline

          &                       & \multicolumn{3}{c}{{\sc cuts}}                                           & \multicolumn{4}{c}{{\sc runtime}}                                                                    & \textbf{}               & \textbf{}              & \textbf{}                          \\
name  & \multicolumn{1}{c}{$\vert V \vert $} & \multicolumn{1}{c}{trian.} & \multicolumn{1}{c}{LK-H} & TSP & trian. & \multicolumn{1}{c}{LK-H} & \multicolumn{1}{c}{TSP} & \multicolumn{1}{c}{total} & \multicolumn{1}{c}{BC n.} & \multicolumn{1}{c}{LB} & \multicolumn{1}{c}{UB}             \\
\noalign{\smallskip}\hline

gr24      & 24                    & 866                        & 25                          & 1       & 0.00   & 0.5                         & 60.6                        & 68.2                      & 936                     & 820.0                  & 820.0                              \\
bayg29    & 29                    & 591                        & 15                          & 1       & 0.00   & 1.1                         & 505.5                       & 508.0                     & 281                     & 843.0                  & 843.0                              \\
bays29    & 29                    & 930                        & 33                          & 1       & 0.00   & 1.3                         & 238.9                       & 242.1                     & 740                     & 814.0                  & 814.0                              \\
dantzig42 & 42                    & 1750                       & 82                          & 3       & 0.00   & 11.8                        & 542.9                       & 586.0                     & 2801                    & 869.5                  & 869.5                              \\
swiss42   & 42                    & 2997                       & 25                          & 1       & 0.00   & 7.9                         & 9683.0                      & 9723.2                    & 1643                    & 877.0                  & 877.0                              \\
gr48      & 48                    & 53310                      & 431                         & 4       & 0.00   & 97.9                        & 10369.5                     & 20085.4                   & 110488                  & 818.5                  & 818.5                              \\
hk48      & 48                    & 10310                      & 329                         & 3       & 0.00   & 89.4                        & 11.9                        & 678.6                     & 14469                   & 823.0                  & 823.0                              \\
eil51     & 51                    & 7709                       & 163                         & 4       & 0.00   & 40.0                        & 14230.2                     & 14851.8                   & 8301                    & 778.5                  & 778.5                              \\
brazil58  & 58                    & 81                         & 2                           & 1       & 0.00   & 1.8                         & 8314.9                      & 8324.6                    & 1                       & 858.0                  & 858.0                              \\
st70      & 70                    & 62803                      & 386                         & 3       & 0.01   & 128.9                       & 12740.9                     & 44798.8                   & 110663                  & 779.5                  & 779.5                              \\
pr76      & 76                    & 57598                      & 969                         & 1       & 0.26   & 532.6                       & 9527.8                      & 86000.0                   & 172533                  & \textit{778.5}         & \multicolumn{1}{r}{\textit{780.0}} \\
eil76      & 76 & 107629 & 631 & 0 & 0.18 & 351.7 & - & 86000.0 & 255729 & \textit{780.0}  & \multicolumn{1}{r}{\textit{-}}  \\
rat99     & 99                    & 33369                      & 30                          & 0       & 3.07   & 37.4                        & -                         & 86000.0                   & 28542                   & \textit{814.5}         & \multicolumn{1}{r}{\textit{-}}     \\
kroB100   & 100                   & 65429                      & 62                          & 0       & 2.81   & 72.7                        & -                         & 86000.0                   & 38288                   & \textit{769.0}         & \multicolumn{1}{r}{\textit{-}}     \\
kroC100   & 100                   & 93370                      & 26                          & 0       & 2.02   & 24.3                        & -                         & 86000.0                   & 40141                   & \textit{775.5}         & \multicolumn{1}{r}{\textit{-}}    \\
\noalign{\smallskip}\hline
\end{tabular}
\end{adjustbox}
\end{table}

\subsection{Generating hard instances by sampling}\label{sec:sampled}
We have also generated a collection of instances for the {\tt Hard-TSPLIB} by using the random sampling procedure discussed in Section \ref{sec:sampling}.
First, we run Algorithm \ref{alg:herSEP} to generate a random set $R$ of vertices of $P_{SEP}$, for a fixed size $n$ of the TSP.
When generating a random vertex $\bm{\bar{x}}^{(h)}$ in Algorithm \ref{alg:herSEP}, we also store the cost vector $\bm c_0^{(h)}$ sampled from the metric cone which yields the vertex $\bm{\bar{x}}^{(h)}$.
Hence, we can compute the integrality gap of the initial {\it easy} TSP instances.
Later, for each random vertex $\bm{\bar{x}}^{(h)}$ in $R$, we solve the IH-OPT$(\bm{\bar{x}}^{(h)})$ problem to get an instance with a larger integrality gap.

Table \ref{tab:hardSampling} reports the three hardest instances we are able to provide for each $n \in \{10, 15, 20,\\25, 30, 35, 40\}$, generating 10 random vertices for each value of $n$.
The table reports first the integrality gap of the sampled cost vector $\bm c_0$, of the optimal cost vector $\bm c^*$ obtained after solving IH-OPT$(\bm{\bar{x}}^{(h)})$, and the gap $\alpha^*_n$ conjectured in \cite{benoit2008finding}, that is, to the best of our knowledge, the highest integrality gap available in the literature.
Then, in the remaining columns, the table reports the average runtime (with the standard deviation) and the average number of BC nodes for solving the instances to optimality using Concorde, first for the instance $\bm c_0$ (random cost vector) and then for $\bm c^*$ (the optimal solution of IH-OPT).
Similarly to the results of the {\tt Hard-TSPLIB}, Table \ref{tab:hardSampling} shows that we are able to generate very hard instances by solving the IH-OPT problem.

\begin{figure}
    \centering
    \includegraphics[scale=0.5]{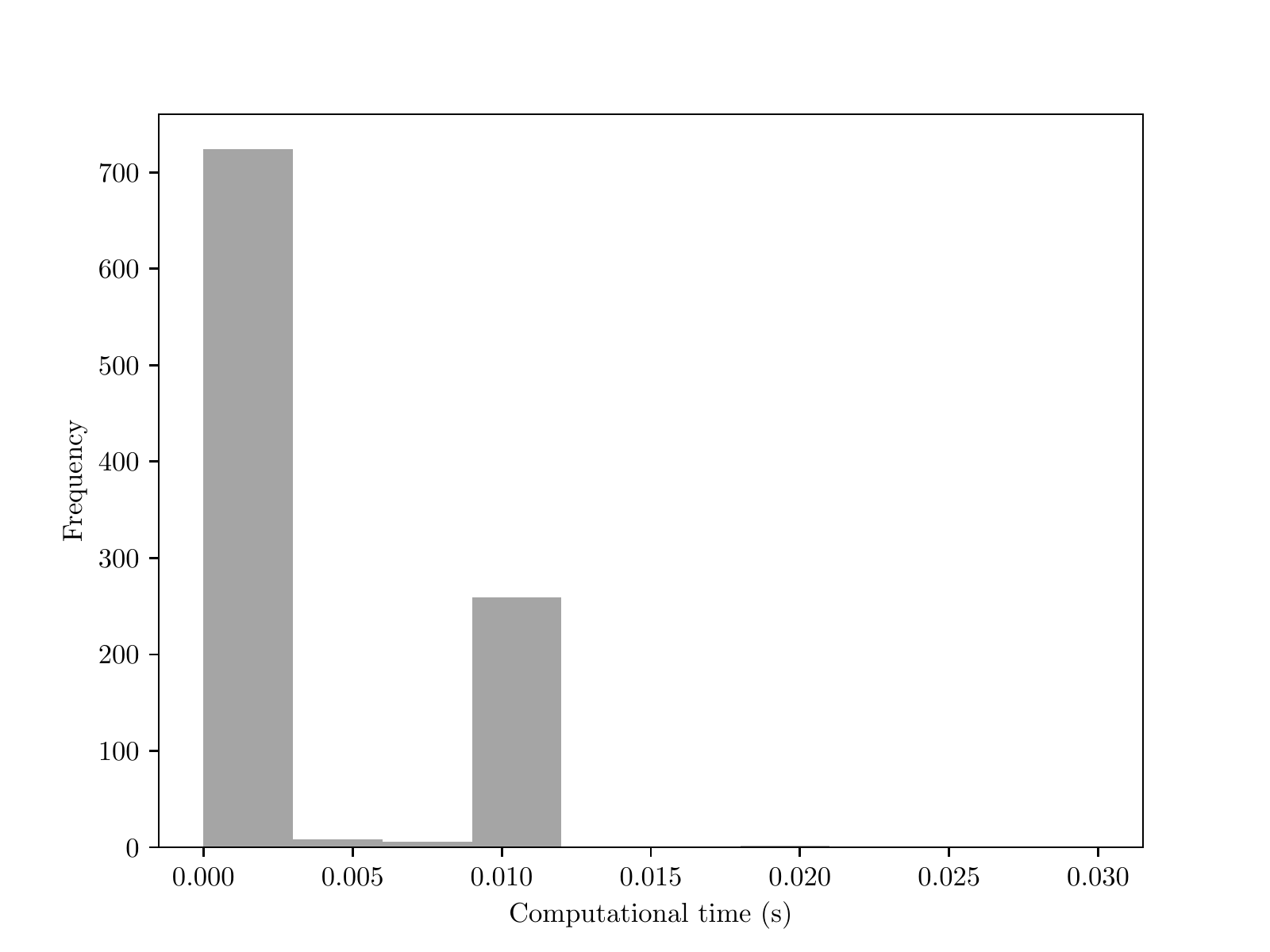}
    \includegraphics[scale=0.5]{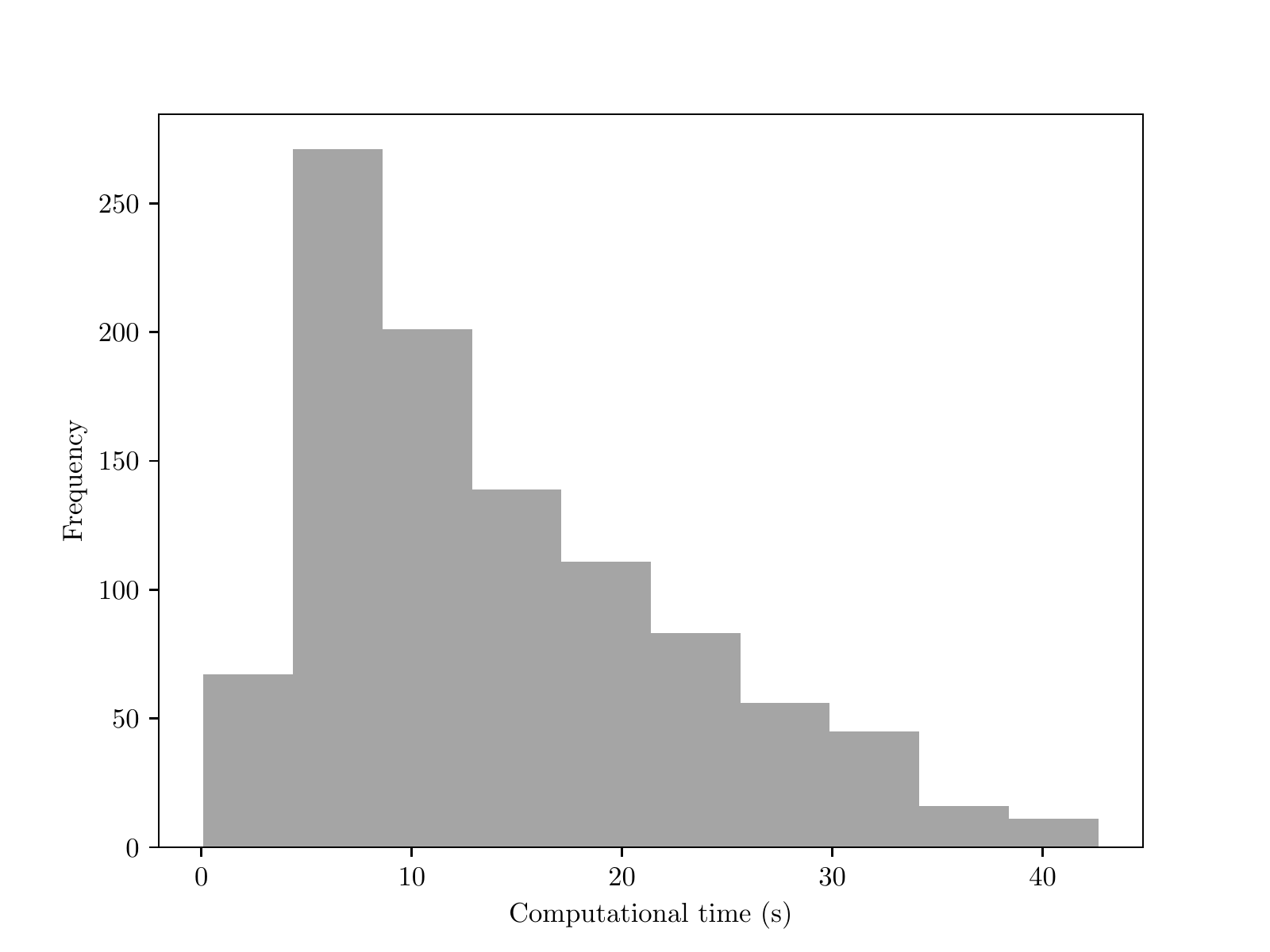}
    \caption{Top: distribution of the computational times for $n = 20$ of 1000 sampled TSP instances. Bottom: distribution of the computational times for $n = 20$ of 1000 instances obtained by the IH-OPT procedure. }
    \label{fig:hist}
\end{figure}

We run a second experiment to study the effect of generating 1000 random instances for $n=20$.
For each random instance $\bm c_0^{(h)}$, with $h = 1,\ldots, 1000$,  we measure the average solution runtime of Concorde.
Then, using the corresponding vertex $\bm{\bar{x}}^{(h)}$ and solving IH-OPT$(\bm{\bar{x}}^{(h)})$, we generate the harder instance, and we measure the runtime again.
Figure \ref{fig:hist}(a) shows the runtime distribution for the random instances $\bm c_0^{(h)}$, while Figure \ref{fig:hist}(b) shows the runtime distribution for $\bm c^*$.
In the top plot, we have a runtime close to zero: we barely get close to 0.030 seconds.
On the other hand, for the bottom plot, we have mean runtime of 10 seconds, with a maximum of around 50 seconds.
In practice, if we sample a large number of vertices we are able to get very challenging small instances. 

\begin{table}[t]
\caption{Results for the {\tt Hard-TSP} instances generated by sampling.}
\label{tab:hardSampling}       
\begin{adjustbox}{width=\textwidth}
\begin{tabular}{crrrrrrrrr}
\hline\noalign{\smallskip}
 & \multicolumn{3}{c}{{ Integrality Gap} - { $\alpha_n(\bm{c})$}} & \multicolumn{3}{c}{{Initial instance} - $\bm c_0$} & \multicolumn{3}{c}{{Hard instance} - $\bm c^*$} \\
$|V|$   & $\bm{c}_0$ & $\bm{c}^*$  & $\alpha^*_n$ \cite{benoit2008finding} &  runtime & stdev & BC n. & runtime & std.dev & BC n. \\
\noalign{\smallskip}\hline\noalign{\smallskip}
10 & 1.007 & 1.153 & 1.176 & 0.00 & 0 & 1 & 0.29    & 0.06    & 5.2     \\
10 & 1.002 & 1.153 & 1.176 & 0.00 & 0 & 1 & 0.23    & 0.06    & 3.8     \\
10 & 1.001 & 1.157 & 1.176 & 0.00 & 0 & 1 & 0.20    & 0.07    & 3.4     \\
15 & 1.005 & 1.172 & 1.222 & 0.00 & 0 & 1 & 4.45    & 1.17    & 12.8    \\
15 & 1.018 & 1.171 & 1.222 & 0.00 & 0 & 1 & 1.68    & 0.48    & 1.2     \\
15 & 1.007 & 1.170 & 1.222 & 0.00 & 0 & 1 & 3.45    & 1.29    & 3.0     \\
20 & 1.007 & 1.179 & 1.246 & 0.01 & 0 & 1 & 30.76   & 7.14    & 66.0    \\
20 & 1.009 & 1.212 & 1.246 & 0.01 & 0 & 1 & 8.86    & 2.15    & 9.2    \\
20 & 1.011 & 1.217 & 1.246 & 0.00 & 0 & 1 & 18.91   & 7.51    & 51.2    \\
25 & 1.022 & 1.244 & 1.262 & 0.01 & 0 & 1 & 140.43  & 32.58   & 383.2   \\
25 & 1.020 & 1.256 & 1.262 & 0.01 & 0 & 1 & 138.66  & 13.31   & 385.2   \\
25 & 1.011 & 1.244 & 1.262 & 0.01 & 0 & 1 & 137.42  & 24.74   & 390.8   \\
30 & 1.002 & 1.002 & 1.273 & 0.01 & 0 & 1 & 195.12  & 195.13   & 578.00  \\
30 & 1.012 & 1.245 & 1.273 & 0.01 & 0 & 1 & 483.10  & 96.62   & 1426.8  \\
30 & 1.000 & 1.263 & 1.273 & 0.01 & 0 & 1 & 156.55  & 21.13   & 454.4   \\
35 & 1.013 & 1.263 & 1.281 & 0.01 & 0 & 1 & 2498.79 & 639.69  & 6565.0  \\
35 & 1.018 & 1.276 & 1.281 & 0.02 & 0 & 1 & 1804.68 & 391.62  & 5193.5  \\
35 & 1.012 & 1.244 & 1.281 & 0.01 & 0 & 1 & 1466.32 & 265.40  & 3962.1  \\
40 & 1.002 & 1.210 & 1.287 & 0.01 & 0 & 1 & 224.85  & 102.18  & 514.0  \\
40 & 1.005 & 1.280 & 1.287 & 0.01 & 0 & 1 & 5889.35 & 823.20  & 13212.4 \\
40 & 1.011 & 1.278 & 1.287 & 0.02 & 0 & 1 & 6279.78 & 874.02 & 14261.0 \\
\noalign{\smallskip}\hline
\end{tabular}
\end{adjustbox}
\end{table}

Table \ref{tab:lightopttimes} reports the average runtime (in seconds) to solve the IH-OPT problem on 1000 different vertices sampled for each $n$ from 10 to 20.
The solution to this problem exhibits a very large runtime variability.
For instance, with $n = 20$, the minimum runtime is of 2 seconds, while the maximum is of 3445 seconds, that is, three orders of magnitude larger.
In addition, we have fitted a linear regression on the $\log_{10}$ of the runtime, getting the following equation
\[
    t = 10^{0.146\cdot n -1.86}.
\]
Using this regression, we try to predict the runtime for each $n$.
For example, with $n = 40$, we would expect a runtime of around 3 hours.
However, we have not observed any significant correlation between the runtime for solving the IH-OPT problem and the average runtime for solving the hard instance with Concorde.

\begin{table}
\centering
\caption{Average computational time and standard deviation of 1000 $\text{IH-OPT}$ on different $\bm{x}$ at a fixed $n$ number of nodes.}
\label{tab:lightopttimes}       
\begin{tabular}{lrrrr}
\hline\noalign{\smallskip}
$n$ & avg. time& std. dev &  min runtime &  max runtime\\
\noalign{\smallskip}\hline\noalign{\smallskip}
10 & 0.63 & 0.39 & 0.29 & 6.86 \\
11 & 0.51 & 0.22 & 0.21 & 1.98 \\
12 & 0.74 & 0.33 & 0.31 & 2.84 \\
13 & 1.19 & 0.64 & 0.45 & 9.08 \\
15 & 2.40 & 1.84 & 0.63 & 26.97 \\
16 & 3.38 & 2.40 & 1.00 & 27.03 \\
17 & 5.82 & 6.41 & 1.29 & 100.14 \\
18 & 9.37 & 22.13 & 1.81 & 622.00 \\
19 & 15.63 & 107.89 & 2.15 & 3396.61 \\
20 & 15.63 & 109.43 & 2.24 & 3445.62 \\
\noalign{\smallskip}\hline
\end{tabular}
\end{table}

\subsection{A comparison with the 3D-Rectilinear instances}
To the best of our knowledge, the current hardest instances of the literature are the Rectilinear 3D instances provided by Zhong \cite{zhong2021lower}.
In this work, the author proposes a family of 3-paths instances with a different configuration of nodes on the three paths and suggests the node configuration that makes them hard to solve for Concorde.
For this reason, we tried to generate, using our computational procedure, metric TSP instances as hard as those 3D Rectilinear instances.
We have tried the following two strategies.

\begin{enumerate}[(a)]
    \item We massively sampled the metric polytope for small $n$ with the algorithm described in Section \ref{sec:sampling} and then apply IH-OPT to the correspondent vertices, until we find an instance with an average runtime competitive with the 3D-Rectilinear instances.
    We call this found instance $\bm{c}_n^S$.
    \item We tried to use our computational procedure to generate hard instances starting from the optimal solution of the instances provided by Zhong \cite{zhong2021lower}.
    Let $\bm{c}^R_n$ the cost vector of a 3D rectilinear instance with $n$ nodes \cite{zhong2021lower}, and let $\bm{x}^R_n$ the corresponding optimal solution of SEP$(\bm{c}^R_n)$.
    Similarly, let $\bm{c}^{IH}_n = \arg\min \mbox{IH-OPT} (\bm{x}_n^R)$ and $\bm{x}^{IH}_n = \arg\min \mbox{SEP}(\bm{c}^{IH}_n)$.
\end{enumerate}
We compared the average runtime of solving TSP$(\bm{c}^R_n)$, TSP$(\bm{c}^{IH}_n)$, and TSP$(\bm{c}^S_n)$.
We have observed from our preliminary tests that Concorde has great runtime variance in solving the same instance:
See, for example, the standard deviation values in Table \ref{tab:hardTSPLIB} and Table \ref{tab:hardSampling}, where we performed 10 independent runs on each instance, using every time a different seed.
This variability is typical of Mixed Integer Programming and branch and cut algorithm (see e.g \cite{fischetti2016improving}, \cite{lodi2013performance}).
For this reason, the runtime comparison instance by instance is meaningless unless the difference is of at least one order of magnitude.
Thus, we present the runtime comparison using linear regression on the logarithm ($\log_{10}$) of the runtime measured in seconds (again, using 10 independent runs for each instance).

\begin{figure}[t!]
\centering
  \includegraphics[scale=0.7]{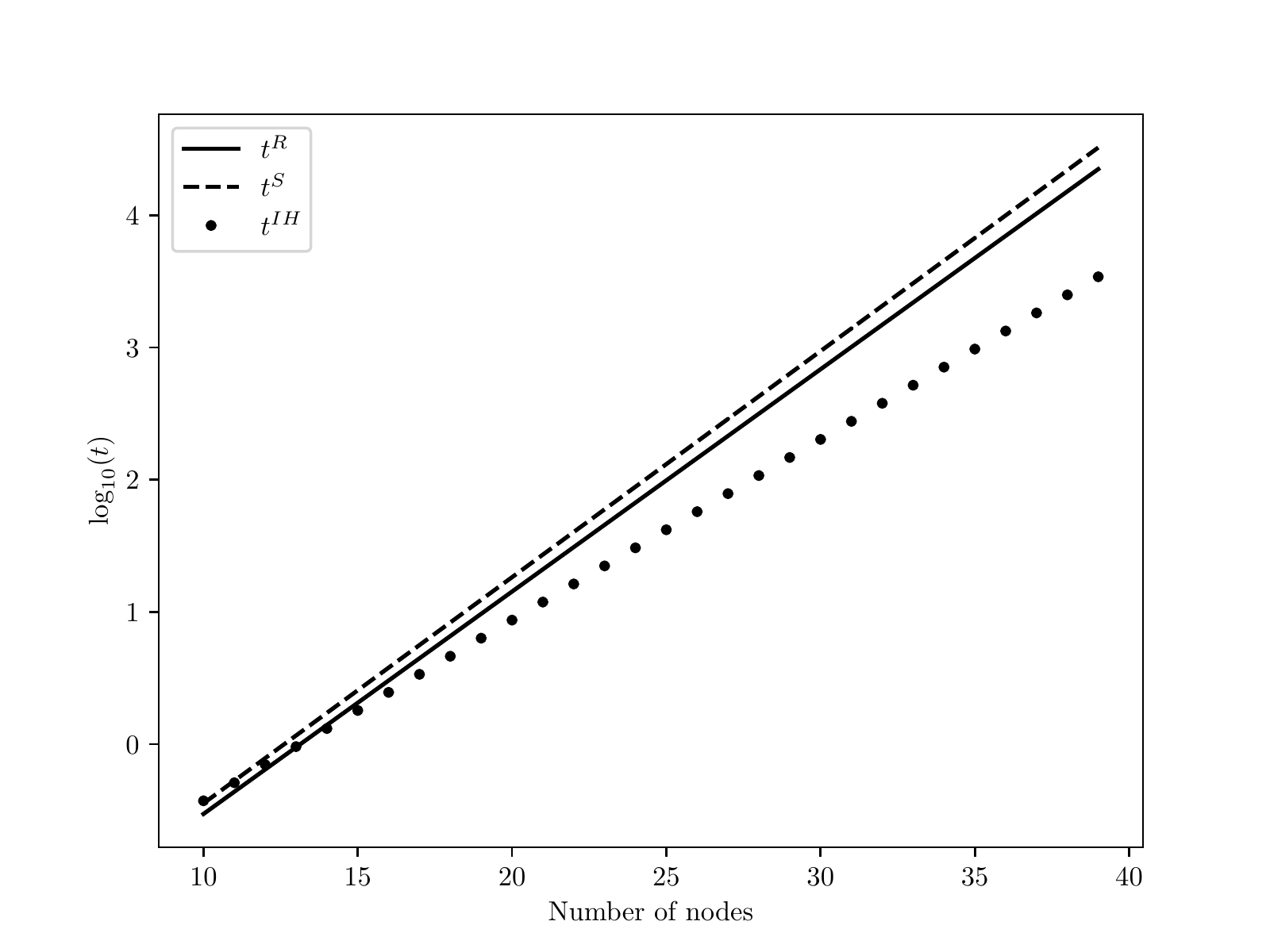}
  \caption{Linear regression for the runtime $t^R_n$, $t^S_n$, and $t^{IH}_n$, corresponding to 3D rectilinear instances, sampled instances, and instances which are solution of IH-OPT$(\bm{x}_n^R)$, respectively.\label{fig:regrZ}}
\end{figure}

Figure \ref{fig:regrZ} shows the three regression lines for comparing the runtime for solving (i) the 3D rectilinear instances, marked with R, (ii) the instances produced using the sampling procedure mentioned above in point (a), marked with S, and (iii) the instances used with the procedure described in the point (b), marked with IH.
The three runtimes are denoted by $t^{R}_n$, $t_n^S$, $t^{IH}_n$, respectively.
The fitted regression lines are the following:
\[
  t^{R}_n = 10^{ 0.168\, n -2.209}, \quad
  t^S_n = 10^{ 0.171\, n -2.157  }, \quad
  t^{IH}_n = 10^{ 0.137\, n -1.794}.
\]
By looking at the regression, the instances that we have found by massively sampling the space are a bit harder than the Rectilinear 3-D instance.
On the contrary, the instances directly obtained from vertices are easier.
Interestingly, in the logarithm scale, the three instances have all the same order of magnitude:
by looking at Figure \ref{fig:regrZ} it is possible to observe that the three families are competitive to each other, with the differences mostly related to the variability of Concorde.
Noteworthy, by observing the structure of costs of the instances created by hand as \cite{zhong2021lower} with respect to the one obtained by the two heuristic procedures, we note less regularity.
This fact will be discussed wider in the next subsection.
In terms of the integrality gap, we do not notice any significant evidence. 
For $n = 12$, the instance we sample has an integrality gap of 1.164 while the Rectilinear 3-D instance has 1.193  and the two runtimes are competitive.
On the contrary, for $n = 20$, we sample an instance with an integrality gap of 1.242, while the Rectilinear 3-D instance has 1.240 and yet the two instances are competitive.
Lastly, we run the computation as reported by Zhong, namely by multiplying the instances by 1000 and rounding to the nearest integer. 
We verified that this procedure might lead to non-metric instances as, due to rounding errors, triangular inequalities might not be satisfied.

\subsection{Structure of small hard instances}\label{sec:study}

The hard instances recently introduced in \cite{hougardy2020hard} and \cite{zhong2021lower} are characterized, by construction, by half-integer solutions for SEP solution, and by support graphs having two triangles where each edge has a weight equal to $\frac{1}{2}$.
Hence, we have looked at the structure of the support graph of our hard instances.
For some instances, such as  {\tt bays29} and {\tt eil51}, the optimal vertex moves from a complicated structure to a 3-path configuration with two triangles having $\frac{1}{2}$-vertices.
For the instance {\tt gr48}, in the original version, the vertices have entries of values in the set $\{0.0, 0.25, 0.5, 0.75, 1.0\}$, while in the corresponding hard version, they have entries in the set $\{0.0, 0.5, 1.0\}$.

\begin{figure}[t!]
    \centering
    \includegraphics[width=0.5\textwidth]{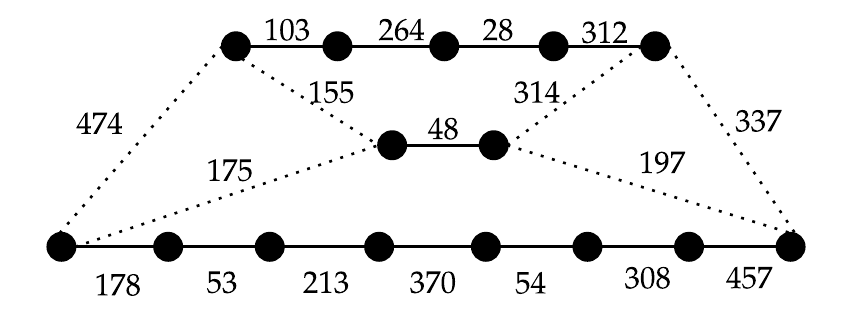}
    \includegraphics[width=0.5\textwidth]{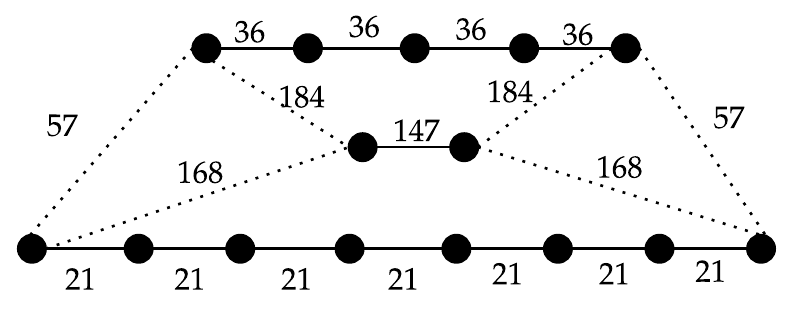}
    \caption{Easy (top) and hard (bottom) support graph associated to one of the sampled vertex with $n = 15$. The cost value on the edge of the support graph is also provided.}
    \label{fig:15007}
\end{figure}

\begin{figure}[t!]
    \centering
    \includegraphics[width=0.7\textwidth]{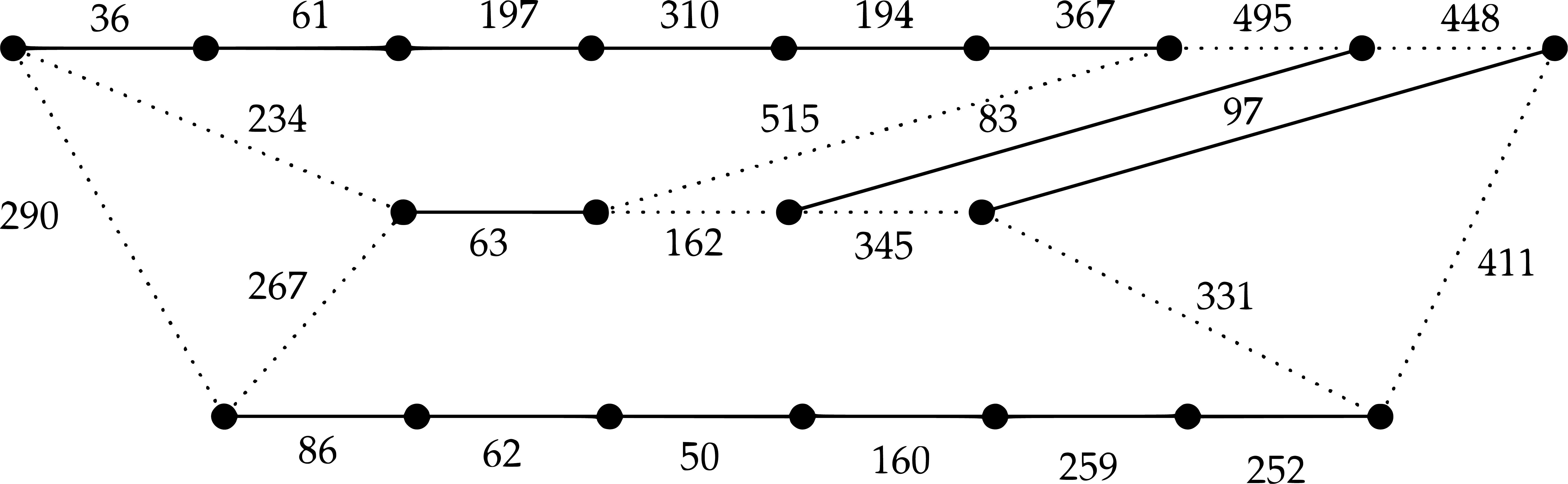}
    \includegraphics[width=0.8\textwidth]{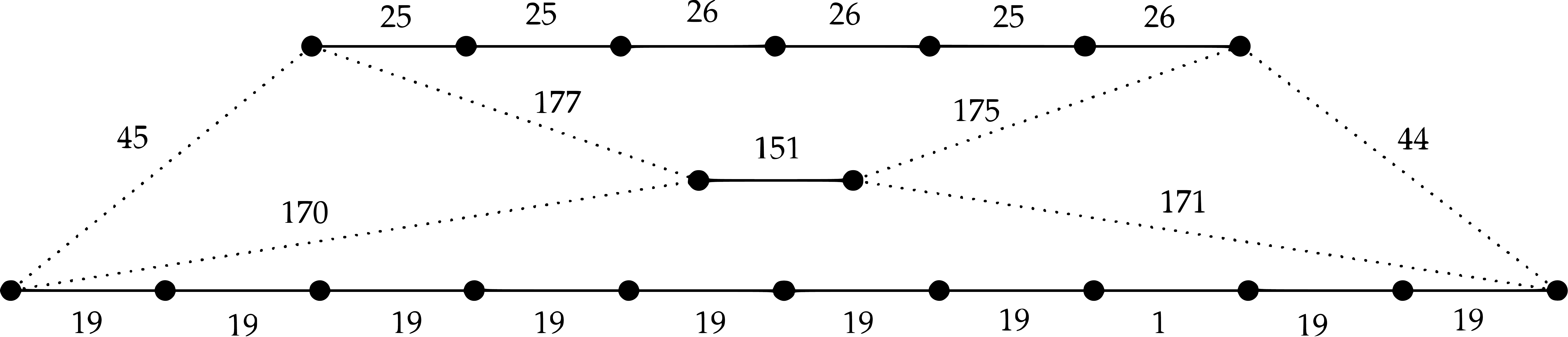}
    \caption{Easy (top) and hard (bottom) support graph associated to one of the sampled vertex with $n = 20$. The cost value on the edge of the support graph is also provided.}
    \label{fig:20012}
\end{figure}

\begin{figure}[t!]
    \centering
    \includegraphics[width=0.8\textwidth]{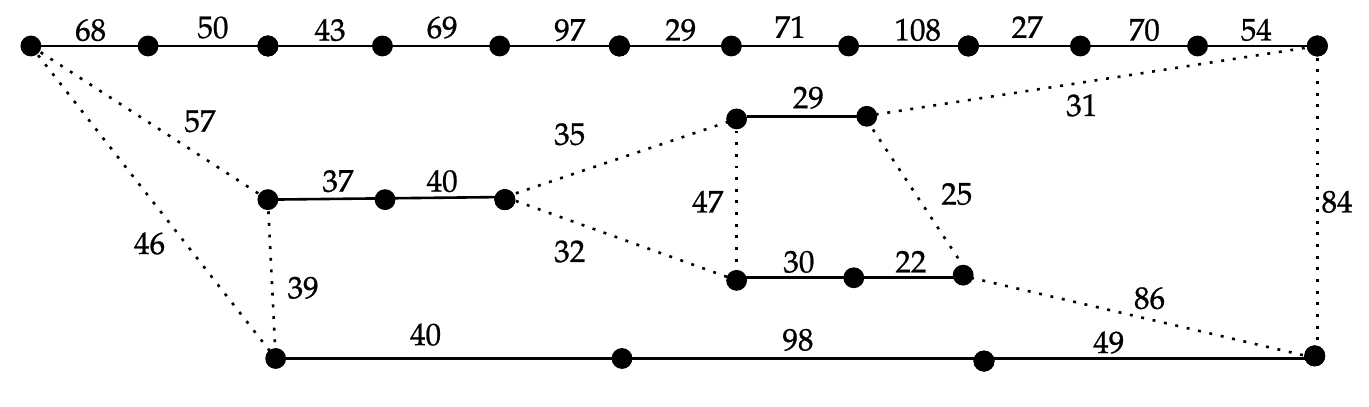}
    \includegraphics[width=0.82\textwidth]{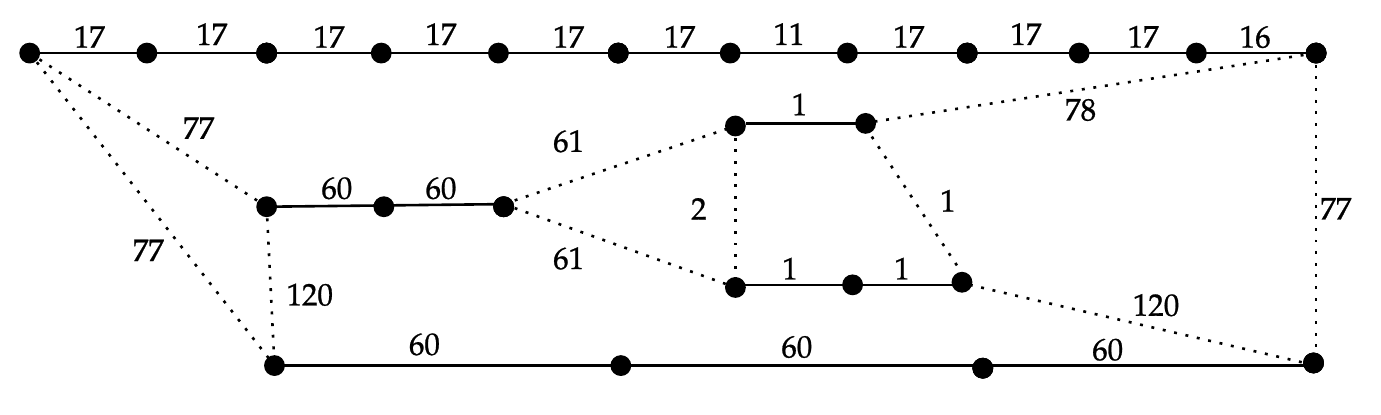}
    \caption{Easy (top) and hard (bottom) support graph associated to the instance gr24 of the TSPLIB. The cost value on the edge of the support graph is also provided.}
    \label{fig:gr24}
\end{figure}

In addition, for three small instances, namely {\tt gr24} and two sampled instances with $n=15$ ($s_{15}$) and $n=20$ ($s_{20}$ ), we studied ``by hand'' the cost structure of the support graph of the optimal solution of SEP.
We selected these 3 instances since they are challenging for Concorde.
The motivation of this study is to investigate if the hard instances share some common cost patterns and/or structures of the support graph.
Figures \ref{fig:15007}, \ref{fig:20012}, and \ref{fig:gr24} show the support graph of the {\it easy} and {\it hard} version of the two instances, where the easy refer to the original cost vector $\bm c_0$ with its SEP vertex solution, and the hard to the optimal solution $\bm c^*$ of IH-OPT.
In the support graphs, the dotted edges correspond to the solution of SEP having value $x_e = \frac{1}{2}$, while the solid edges correspond to $x_e = 1$.
The missing edges have $x_e = 0$.
The label of each edge gives its cost in the TSP instance.\\

We observe that the main difference between the easy and hard instances is in the pattern of the edge length.
In the easy instances, the edge costs look randomly distributed, while in the hard instances obtained after solving IH-OPT there is a clear cost pattern.
The edges lying along the same path have nearly the same cost, while the dashed edges connecting two distinct paths generally have a cost nearly equal to the sum of the costs of a single edge on a path.
In addition, the overall length of the 3-paths are almost equal: for instance, in $s_{15}$ (Figure \ref{fig:15007}), in the easy version, we have on each path a sum of, respectively, 707, 48, 1633, while in the hard one, we find 144, 147, 147.
Finally, notice that since in IH-OPT we have removed the slackness constraints introduced in OPT$_h$, we do not have any guarantee that the optimal solution for the SEP associated to $\bm{c}^0$ and $\bm{c}^*$ is the same.
However, in $s_{15}$ we have the $\bm{x}^0 = \bm{x}^*$, while for $s_{20}$ and {\tt gr24} we get $\bm{x}^* \neq \bm{x}^0$, as shown in Figure \ref{fig:20012} and Figure \ref{fig:gr24}.

\section{Conclusions}
In this work, we have introduced a computational procedure to generate metric TSP instances which have large integrality gaps and are challenging for Concorde, the state-of-the-art TSP solver.
As a by-product, we have introduced the {\tt Hard-TSPLIB}, a collection of small but challenging metric TSP instances, which are not generated explicitly exploiting specific cost structures, as in \cite{hougardy2020hard,zhong2021lower}.
Notice that, to the best of our knowledge, all the hard instances from the literature have half-integer optimal SEP solutions.
On the contrary, the instances of the {\tt Hard-TSPLIB} have a larger variety of fractional optimal vertices (i.e., they are not only half integers).
We expect our new instances will serve as a benchmark for designing new exact and heuristic methods for solving the TSP problem.

Curiously, we have observed that the most challenging instances generated using our computational procedure have regular cost patterns, with several edges sharing the same costs, and several paths on the support graphs having the same length.
These types of cost patterns are in common with the manually-generated hard instances recently introduced in \cite{hougardy2020hard} and \cite{zhong2021lower}.
Hence, we believe that our {\tt Hard-TSPLIB} instances could help further studies in the cost structures of TSP instances.

We emphasize that our framework is really general and can be applied to any combinatorial optimization problem.
There are several combinatorial optimization problems that have metric costs, such as, for instance, the metric Steiner tree \cite{bern1989steiner} or the shortest Euclidean minimum spanning tree \cite{fischetti1994weighted}.
The main difference would be to change the TSP constraints into different inequalities for the corresponding combinatorial objects.

\bibliography{mainarXiv}  

\end{document}